\documentclass[12pt]{article}
\usepackage{amsmath,amsopn,amssymb,amsfonts,ctable}
\usepackage{epsfig}
\usepackage{hyperref}
\usepackage{verbatim}
\usepackage{multirow}
\usepackage{arydshln}
\usepackage{amsthm}
\usepackage[margin=1cm]{caption}

\makeatletter \@addtoreset{equation}{section} \makeatother
\makeatletter \@addtoreset{figure}{section} \makeatother

\def\IC{\mathbb{C}}
\def\IN{\mathbb{N}}
\def\IR{\mathbb{R}}\def\IZ{\mathbb{Z}}\def\IQ{\mathbb{Q}}

\def\tr{{\rm tr}}

\def\be{{\bar e}}

\def\bea{\begin{eqnarray}}
\def\eea{\end{eqnarray}}
\def\be{\begin{equation}}
\def\ee{\end{equation}}

\def\={\;  = \;}
\def\+{\, + \,}

\newtheorem{thm}{Theorem}[section]

\newtheorem{proposition}{Proposition}[section]
\newtheorem{exmp}{Example}[section] 
\newtheorem{claim}{Claim}[section]
\newtheorem{rmk}{Remark}[section]
\newtheorem{conj}[thm]{Conjecture} 

\usepackage{cleveref}
\crefname{section}{ÃÂ¤}{ÃÂ¤ÃÂ¤}
\Crefname{section}{ÃÂ¤}{ÃÂ¤ÃÂ¤}

\begin{document}

\begin{titlepage}
\vfill
\begin{flushright}


\end{flushright}
\vfill
\begin{center}
{\Large\bf Traces of Singular Moduli and  Moonshine for the Thompson Group}

\vskip 1cm

Jeffrey A. Harvey and Brandon C. Rayhaun
\vskip 5mm
{\it Enrico Fermi Institute and Department of Physics \\
University of Chicago \\
5620 Ellis Ave., Chicago Illinois 60637, USA}

\end{center}
\vfill

\begin{abstract}
\noindent

We describe  a relationship between the representation theory of the Thompson sporadic group and a weakly holomorphic modular form of weight one-half that 
appears in work of Borcherds and Zagier on Borcherds products and traces of singular moduli. We conjecture the existence of an infinite dimensional graded module for the Thompson group
and provide evidence for our conjecture by
constructing McKay--Thompson series for each conjugacy class of the Thompson group that coincide with weight one-half modular forms of higher level. 
We also observe a discriminant property in this moonshine for the Thompson group that is closely related to the discriminant property conjectured to exist in Umbral Moonshine.
\end{abstract}

\vfill
\end{titlepage}

\renewcommand{\thefootnote}{\#\arabic{footnote}}
\setcounter{footnote}{0}

\section{Introduction}

There appear to be two distinct types of moonshine phenomena in mathematics:  Monstrous moonshine involving the modular function
\be 
J(\tau)= q^{-1} + 196884 q + 21493760 q^2 + \cdots \,
\ee
and the representation theory of the largest sporadic group, the Monster group  (${\mathbb M}$) \cite{conwaynorton,flm,borchmon}; and 
Umbral Moonshine \cite{um,mum} which involves a set of vector-valued mock modular forms $H^X(\tau)$,
labelled by ADE root systems $X$ with common Coxeter numbers and total rank $24$, which exhibit moonshine for  finite groups $G^X=\mathrm{Aut}(L^X)/W(X)$
where $\mathrm{Aut}(L^X)$ is the automorphism group of the Niemeier lattice determined by $X$ and $W(X)$ is the Weyl group of $X$.
Umbral Moonshine extends and generalizes  Mathieu moonshine connecting  the weight one-half mock modular form
\be
H^{(2)}(\tau) = 2 q^{-1/8} \left(-1 + 45 q + 231 q^2 + 770 q^3 + \cdots \right) \, ,
\ee
the Mathieu group $M_{24}=G^X$ for $X=A_1^{24}$ and
the elliptic genus of $K3$ surfaces first observed in \cite{eot} and studied extensively since then including complete computations of the McKay--Thompson series \cite{cheng,GaberdielI, GaberdielII,EguchiI}.
Monstrous moonshine is now best understood in terms of a central charge $c=24$ holomorphic conformal field theory (CFT) or vertex operator algebra (VOA), see  \cite{flmbook} for details,  and there
are related constructions exhibiting moonshine for the Baby Monster group \cite{hohn} and the Conway group \cite{duncancon} also based on CFT constructions. 
While the existence of the modules predicted by Umbral Moonshine has now been proved \cite{gannon, dgo} and there exists an explicit construction of the modules
for the case $X=E_8^3$ \cite{dunhar}, it seems fair to say that much remains to be done to explicate the relationship between Umbral Moonshine and CFT/VOA structures.

There are hints of connections between these two types of moonshine.  For example the ADE root system $X$ was used to define
a set of weight zero modular functions in \cite{mum} and these are also hauptmoduls that appear as McKay--Thompson series
in Monstrous Moonshine, thus giving a correspondence between cases of Umbral Moonshine and conjugacy classes of the Monster group.
Further hints in this direction come from recent work relating dimensions of representations appearing in Umbral Moonshine when $X$
is a pure $A$-type root system  to values of the McKay--Thompson series of Monstrous Moonshine \cite{ono}. 

One of the elements present in Umbral Moonshine that seems to have no analog in Monstrous Moonshine is the discriminant property,
which relates the number fields on which the irreducible representations of $G^X$ are defined, and the discriminants of the
vector-valued mock modular form $H^X$  \cite{um,mum}.  

In this paper we provide evidence for a moonshine phenomenon that shares certain features with both Monstrous moonshine and Umbral Moonshine, 
including a discriminant property similar to one observed in
Umbral Moonshine, but involving 1) modular rather than mock modular forms, and 2) the Thompson sporadic group, a group much larger than
any of the groups of Umbral Moonshine and one with a natural connection to the Monster. 

To introduce the elements playing a role in our analysis we recall some of the results from Zagier's work \cite{zagier} on 
weakly holomorphic modular forms, Borcherds products and traces of singular moduli. Define the Jacobi theta function
\be
\theta(\tau)= \sum_{n \in \IZ} q^{n^2}  \, ,
\ee
with $q=e^{2 \pi i \tau}$, 
which is holomorphic in the upper half plane $\mathfrak{h}$, lives in the ``Kohnen plus-space" of functions with
a Fourier expansion at infinity of the form $\sum c(n) q^n$ with $c(n)=0$ unless $n \equiv 0,1$ modulo $4$, and transforms with a well known multiplier system under the congruence subgroup $\Gamma_0(4)$. Let $M_{1/2}^!$ be the space of functions transforming like
$\theta$ under $\Gamma_0(4)$ and also live in the Kohnen plus-space, but are allowed to be meromorphic at the cusps. In
\cite{zagier} a special basis of $M_{1/2}^!$ is constructed consisting of functions $f_d$ with $d \equiv 0,3$ modulo $4$ with
Fourier expansions of the form
\be \label{fdeqn}
f_d(\tau) = q^{-d} + \sum_{n>0} A(n,d) q^n \, ,
\ee
and with $f_0(\tau)= \theta(\tau)$. The function $f_3$ can be constructed ``by hand" as
\be \label{ftdef}
f_3(\tau)= - \frac{1}{20} \left( \frac{[\theta(\tau), E_{10}(4 \tau)]_1}{\Delta(4 \tau)} + 608~ \theta(\tau) \right) \, 
\ee
where $E_{10}$ is the weight ten Eisenstein series and $\Delta$ is the weight twelve cusp form. 
Equation (\ref{ftdef}) involves the first Rankin-Cohen bracket which is defined
for two modular forms $f,g$ of weight $k,l$ respectively as
\be
[f,g]_1= k  fDg-l g Df
\ee
with $D= \frac{1}{2 \pi i} \frac{d}{d \tau} $. An effective algorithm for computing the remaining $f_d$ is given in \cite{zagier}, but
will not be needed in what follows. Explicit computation  gives the Fourier expansion
\begin{align} 
f_3(\tau) &= q^{-3}  -248 q + 26752 q^4 - 85995 q^5 + 1707264 q^8 -  \\
& \ \ \ \ \ \ \ \ \ \ \ \ 4096248 q^9 +44330496 q^{12} -91951146q^{13} + 708938752q^{16} + \cdots  \, .\nonumber
\end{align}
The coefficients of $f_3$ appear to have a close connection to dimensions of irreducible representations of
the sporadic Thompson group Th of order  $2^{15} \cdot 3^{10} \cdot 5^3 \cdot 7^2 \cdot 13 \cdot 19 \cdot 31$ $ \sim 9  \cdot 10^{16}$. 
For example, $248$, $85995$, $1707264$ and $44330496$ are dimensions of Th irreducible representations while $26752=27000-248$ and $4096248=4096000+248$
with $27000$ and $4096000$ also dimensions of Th irreducible representations.

There is a well understood connection between $f_3$ and the Thompson group arising from the
following facts.  According to the generalized moonshine conjecture of Norton \cite{norton}, for each pair of commuting elements $(g,h)$ in ${\mathbb M}$
there exists a generalized character $Z(g,h,\tau)$  such that $Z(1,h,\tau)$ are the McKay--Thompson series $T_{[h]}$ of \cite{conwaynorton}, 
\be
Z(g^a h^c, g^b h^d, \tau)= \psi Z(g,h,\tfrac{a \tau+b}{c \tau+d})
\ee
for some constant $\psi$ with 
$\begin{pmatrix} a & b \\ c & d \end{pmatrix} \in \mathrm{SL}(2,\IZ)$,
and such that the coefficients of the $q$-expansion of $Z(g,h,\tau)$ for fixed $g$ form characters
of a graded representation of a central extension of the centralizer of $g$ in ${\mathbb M}$. According to \cite{dgh} the characters $Z(g,h,\tau)$ have an interpretation as traces over representations
of the centralizers of $g$ acting in a Hilbert space twisted by $g$, that is a twisted module  of the Monster vertex operator algebra $V^\natural$.

In particular, the character $Z(g,1,\tau)$ for $g$ in the $3C$ class of ${\mathbb M}$ 
is the Borcherds lift \cite{borcherds} of $f_3$:
\bea
Z(3C,1,\tau) = j(\tau)^{1/3} &= & q^{-1/3}  \prod_{n>0} (1-q^{n})^{A(n^2,3)} \\
 &=& q^{-1/3} + 248 q^{2/3} +4124 q^{5/3} + 34752 q^{8/3} + \cdots \, .
\eea
As a result the Thompson group
has a natural action on the $3C$-twisted module of the Monster module $V^\natural$, see \cite{carnahan} for further details.  Note  however that the lift
involves only the coefficients $A(n,3)$ with $n$ a perfect square, and thus this connection does not provide an explanation for why the
coefficients of non-square powers of $q$ in $f_3$ also exhibit a connection with the Thompson group.

The Thompson moonshine structure, as well as a closer similarity to elements of Umbral moonshine, is brought out by considering
instead of $f_3$ the function
\be
{\cal F}_3(\tau) = 2 f_3(\tau) + 248 \theta(\tau) =  \sum_{\substack{m=-3 \\ m \equiv 0,1 ~{\rm mod} 4}}^{\infty}  c(m) q^m 
\ee
which is also an element of $M_{1/2}^!$. The coefficients $c(m)$ are given in Table 1 for $-3 \le m \le 33$.
We then observe using the character table for the Thompson group provided in Tables \ref{thchartabone}--\ref{thchartabfour} that  each coefficient $c(m)$ for $m \le 12$ can be
interpreted, up to sign,  as the dimension of  either a single real irreducible representation (with multiplicity) or a representation of the form $V\oplus \overline V$ where
$V$ is an irreducible representation and $\overline V$ is the conjugate representation.  The coefficients of higher powers of $q$ have more complicated decompositions into
irreducible representations with positive integer coefficients which will be determined later and can be found tabulated in Tables \ref{decomp1} and  \ref{decomp2}.

\begin{table}\begin{center}
\begin{tabular}{ccc}\toprule
$k$ & $c(4k)$ & $c(4k+1)$ \\\midrule  
-$1$ & $0$ & $2 \cdot 1$ \\
$~~0$ & $248$ & $0$  \\
$~~1$ & $2 \cdot 27000$ & -$2 \cdot  85995$  \\
$~~2$ & $ 2 \cdot 1707264$ & -$2 \cdot 4096000$  \\
$~~3$ & $2 \cdot 44330496$ & -$2 \cdot 91951146$  \\
$~~4$ & $2 \cdot 708939000 $ & -$2 \cdot 1343913984$  \\
$~~5$ & $2 \cdot 8277534720$ & -$2 \cdot 14733025125$  \\
$~~6$ & $2 \cdot77092288000 $ & -$2 \cdot 130880765952 $  \\
$~~7$ & $2 \cdot 604139268096$ & -$2 \cdot  988226335125 $ \\
$~~8$ & $2 \cdot 4125992712192$ & -$2 \cdot  6548115718144$ \\
\bottomrule
\end{tabular}
\end{center}\caption{\label{thomqi} {Coefficients in the Fourier expansion of ${\cal F}_3$ .}}
\end{table}

To be more specific, label the irreducible representations of Th as $^{d_i}V_i$, $i=1,\dots, 48$ with $d_i$ the dimension of the irreducible representation, a list of which can be inferred from the character table of the Thompson group provided in Tables \ref{thchartabone}--\ref{thchartabfour}, and abbreviate this at times to $V_i$. The Fourier coefficients of $\mathcal{F}_3$ up to $q^{12}$  then imply the relationship to the irreducible representations of Th summarized in Table \ref{thomqk}. 

\begin{table}\begin{center}
\begin{tabular}{rc}\toprule
$c(k)$ & $ {\rm Decomposition} $  \\\midrule  
$c(\text{-} 3)$ & $ 2 \cdot ^{1}V_1 $ \\
$c(0)$ & $ ^{248}V_2$ \\
$c(4)$ & $  ^{27000}V_4 \oplus  {^{27000} V_5} $ \\
-$c(5)$ & $^{85995}V_9 \oplus  {^{85995}V_{10}} $ \\
$c(8)$ & $^{1707264}V_{17} \oplus  {^{1707264}V_{18}} $ \\
-$c(9)$ & $ ^{4096000}V_{22}  \oplus  {^{4096000}V_{23}} $ \\
$c(12)$ & $ 2 \cdot {^{44330496}V_{40}} $ \\
\bottomrule
\end{tabular}
\end{center}\caption{\label{thomqk} {Connection between the coefficients $c(k)$ of ${\cal F}_3$ and Thompson representations for $k \le 12$.}}
\end{table}

As usual in the study of moonshine,  this observation suggests that there exists an infinite dimensional $\IZ$-graded module for the Thompson group
\be \label{decomp}
W= \bigoplus_{\substack{m \ge -3 \\ m \equiv 0,1~ {\rm mod} ~4}}^\infty W_m
\ee
where we demand that the module be compatible with the Fourier coefficients of ${\cal F}_3$ in the sense that $|c(m)|= {\rm dim} W_m$ . 

It is natural to associate the alternating signs exhibited by the coefficients $c(m)$ to
a superspace structure. A superspace is a $\IZ/2 \IZ$ graded vector space $V=V^{(0)} \oplus V^{(1)}$ with $V^{(0)}, V^{(1)}$ the even and odd elements
of $V$ respectively. The supertrace of a linear operator $L$ on $V$ that preserves the $\IZ/2 \IZ$ grading is then defined to be
$\mathrm{str}_V L= \mathrm{tr}_{V^{(0)}} L -\mathrm{tr}_{V^{(1)}} L$. If we take $W_m$ for $m \ge 0$ and  $m=0 {\mod 4}$ to have vanishing odd part,
$W_m$ for $m \ge 0$ and $m=1 ~{\mathrm{mod}}~ 4$ to have vanishing even part, and $W_{-3}$ to have vanishing odd part then we demand that $c(m)= \mathrm{str}_{W_m} 1$
for all $m$. Note that the ``wrong" sign of the coefficient of the singular $q^{-3}$ term is similar to the structure exhibited by the mock modular forms
$H^X_r$ of Umbral moonshine which have a singular term $-2 q^{-1/2m}$ with a negative coefficient in the $r=1$ component of $H^X_r$ with $m$ the Coxeter number of $X$
while all the other coefficients are positive integers. 

It is natural to also consider the decomposition of each component of the supermodule $W_m=W_m^{(0)} \oplus W_m^{(1)}$
into irreducible representations of the Thompson group,
\bea \label{thdecom}
W_m^{(a)} &= \bigoplus_{i=1}^{48}  b^{(a)}_{m,i} \cdot {V_i} 
\eea
with $b^{(a)}_{m,i} \in \IN_0$ for $a=0,1$.  The Fourier coefficients of $\mathcal{F}_3$ alone are not enough to determine this decomposition; for instance, on the basis of dimension, one can trivially decompose each $W_m^{(a)}$ into $\dim W_m^{(a)}$ copies of the trivial representation. However, the fact that the $c(m)$ for small values of $m$ can be written so nicely in terms
of the dimensions of non-trivial irreducible representations of Th suggests that the $W_m^{(a)}$ are non-trivial and that there is some structure associated
to the choice of representations.

To test this idea we follow a procedure that is now standard and consider McKay--Thompson series for each conjugacy class $[g]$ in the Thompson group by replacing $c(m)= \mathrm{str}_{W_m} 1 $ in the Fourier
development of ${\cal F}_3$ by $ \mathrm{str}_{W_m}(g)$ where $g$ is any representative of $[g]$.
We thus define the $48$ McKay--Thompson series corresponding to each conjugacy class $[g]$ by
\be
{\cal F}_{3,[g]}(\tau)=  \sum_{\substack{m=-3 \\ m \equiv 0,1 ~{\rm mod}~ 4}}^{\infty} \mathrm{str}_{W_m}(g) q^m \, .
\ee
For an arbitrary choice of decomposition in eqn.(\ref{decomp}) there is no reason to expect that the ${\cal F}_{3,[g]}$ exhibit interesting modular
properties. Conversely, if the ${\cal F}_{3,[g]}$ do exhibit interesting modular properties this should be regarded as evidence for an interesting
relation between the Thompson group and a class of weakly holomorphic weight one-half modular forms. In Monstrous and Umbral moonshine
the analogous McKay--Thompson series are modular forms for the congruence subgroups $\Gamma_0(o(g))$ with $o(g)$ the order of the associated
group element, often with a non-trivial multiplier system.
Thus to test our proposed moonshine connection between ${\cal F}_3$ and the Thompson group we need to generalize ${\cal F}_3$ to 
weakly holomorphic modular forms at level $N$ for $N$ that are orders of Th or multiples of orders of Th. 

The outline of the rest of the paper is as follows. In the second section we recall results relating traces of singular moduli to coefficients of weakly holomorphic modular forms
and more generally to the coefficients of Maass-Poincar\'e series and Rademacher series, both at level $1$ and at level $N$  following \cite{zagier, bringmann, mp, ChengDuncanII}. The third section provides evidence
for the modularity of the McKay--Thompson series by comparing their $q$-expansions to those of weakly holomorphic weight one-half modular forms at level $N$.  We construct the required modular forms
using both traces of singular moduli and coefficients of Rademacher series twisted by multiplier systems similar to those appearing in \cite{um,mum}.  In the fourth section we discuss a discriminant property that relates the discriminant of the
quadratic forms that appear in the computation of the coefficients of ${\cal F}_3$ in terms of traces of singular moduli to the fields over which the representations of Th attached
to these coefficients are defined. We compare and contrast this discriminant property with the discriminant property observed in Umbral moonshine in \cite{um,mum}. An appendix
deals with details of the expressions for coefficients of the weight one-half weakly holomorphic forms appearing here in terms of the computation of traces
of singular moduli at level $N$.

\section{Traces of singular moduli and coefficients of Poincar\'e series}

Here we recall results of Zagier \cite{zagier} that express the coefficients $c(m)$ as traces of singular moduli and results of Bringmann
and Ono \cite{bringmann} derived from a study of Maass-Poincar\'e series that give explicit expressions for the $c(m)$ in terms
of Kloosterman sums. We then extend these results to level $N$ following \cite{mp, ChengDuncanII}. 

\subsection{Notation and preliminaries}
In what follows we often use the notation $e(x)=e^{2 \pi i x}$ and we write $q=e(\tau)$ with $\tau$ in the upper half plane, $\tau \in \mathfrak{h}$.
We define the group
$\Gamma_0(n)$ to be the set of elements 
\be
\begin{pmatrix} a & b \\ c & d  \end{pmatrix} \in \mathrm{SL}(2, \IZ) \text{ with } c \equiv 0 ~{\rm mod}~n
\ee
and also let $\left(\frac{m}{n}\right)$ denote the Kronecker symbol.

For $w\in \IR$ and $\Gamma$ a subgroup of $\mathrm{SL}(2,\IR)$ containing $\pm I$ and commensurable with $\mathrm{SL}(2,\IZ)$ we call a function $\psi:\Gamma\to \IC$ a multiplier system for $\Gamma$ with weight $w$ if 
\be
\psi(\gamma_1)\psi(\gamma_2)\mathrm{j}(\gamma_1,\gamma_2\tau)^{w/2}\mathrm{j}(\gamma_2,\tau)^{w/2} = \psi(\gamma_1\gamma_2)\mathrm{j}(\gamma_1\gamma_2,\tau)^{w/2}
\ee
for each $\gamma_1,\gamma_2\in \Gamma$, where $\mathrm{j}(\gamma,\tau) = (c\tau+d)^{-2}$. We will deal almost exclusively with multiplier systems which depend only on the bottom row of matrices in $\Gamma$ and will thus take the liberty of abusing notation slightly by setting 
\be
\psi(c,d) \equiv \psi\left( \begin{array}{cc} \ast & \ast \\ c & d \end{array}   \right).
\ee
There is the standard $(\psi,w)$-action of $\Gamma$ on holomorphic functions $f:
\mathfrak{h}\to\IC$ on the upper-half plane given by 
\be
\left( f\vert_{\psi,w}\gamma \right)(\tau) \equiv f(\gamma\tau)\psi(\gamma)\mathrm{j}(\gamma,\tau)^{w/2}
\ee
which allows us to define a weakly holomorphic modular form of weight $w$ and multiplier system $\psi$ over $\Gamma$ as a function $f$ which is invariant under this action and holomorphic in the interior of $\mathfrak{h}$, but is allowed to be meromorphic at the cusps. We will typically restrict ourselves to $w=1/2$, $\Gamma = \Gamma_0(4N)$ and refer to the multiplier system associated with $\theta$ as
\be
\psi_0(\gamma) = \left(\frac{c}{d}\right)\epsilon_d
\ee
where
\be
\epsilon_d \equiv 
 \begin{cases}
 1, &  d=1 ~{\rm mod}~4 \\
 i, &  d=3 ~{\rm mod}~4 .
 \end{cases}
\ee

We use the symbol $\chi$ both for genus characters of quadratic forms and for the characters of the Thompson group, but the context should clear up any possible ambiguity. 

\subsection{Weakly holomorphic weight one-half modular forms at level one}

First we describe the relation of the coefficients $A(n,3)$  appearing in equation (\ref{fdeqn}) for $d=3$ to traces of singular moduli following \cite{zagier}.
Let
\be
J(\tau)= q^{-1} + 196884 q + 21493760 q^2 + \cdots
\ee
be the normalized hauptmodul for $\mathrm{SL}(2,\IZ)$, let $d$ be a positive  integer with $d \equiv 0,3 ~{\rm mod}~4$ and denote
by ${\cal Q}_d$ the set of positive definite binary quadratic forms
\be
Q(X,Y)=[a,b,c]=a X^2+b XY + c Y^2, \qquad a,b,c \in \IZ
\ee
with discriminant $b^2-4ac=-d$.  We can define an action of the full modular group $\mathrm{SL}(2,\IZ)$ on quadratic forms in the usual way
\be
Q|_\gamma(X,Y)= Q(pX+qY,rX+sY)
\ee
for $\gamma=\begin{pmatrix} p & q \\ r & s \end{pmatrix} \in \mathrm{SL}(2,\IZ)$ and call two quadratic forms $\Gamma$-equivalent if they are equivalent under the group action $\Gamma$ inherits from the full modular group. 

Since $d>0$, each $Q \in {\cal Q}_d$ has a unique root $\alpha_Q \in \mathfrak{h}$. The value of $J(\alpha_Q)$
depends only on the $\mathrm{SL}(2,\IZ)$ equivalence class of $Q$. The modular trace function ${\bf t}(d)$ is defined as the sum of $J(\alpha_Q)$
over $\mathrm{SL}(2,\IZ)$-equivalence classes weighted by a factor $w_Q$ which is the order of the stabilizer of $Q$ in $\mathrm{PSL}(2,\IZ)$ and is $3$ if $Q$ is $\mathrm{SL}(2,\IZ)$-equivalent to $[a,a,a]$, $2$ if $Q$
is $\mathrm{SL}(2,\IZ)$-equivalent to $[a,0,a]$ and $1$ otherwise:
\be
{\bf t}(d)= \sum_{Q \in {\cal Q}_d /\Gamma} \frac{1}{w_Q} J(\alpha_Q) \, .
\ee
Theorem 1 of \cite{zagier} gives $A(1,3)={\bf t}(3)= -248$ which is, up to a sign, the dimension of an irreducible representation of Th. 

The coefficients $A(n,3)$ when $n>1$ is a fundamental discriminant  are also given by a modular trace function that involves quadratic forms of discriminant $-3n$.  One considers the trace twisted by a genus character $\chi_{n,-3}$ which
assigns to a quadratic form $Q$ of discriminant $-3n$ the value $\pm 1$ determined by 
\be
\chi_{n,-3}(Q)= \left( \frac{n}{p} \right)
\ee
where $p$ is any prime represented by $Q$ and not dividing $3n$.  If in addition $n$ and $3$ are coprime then
\be \label{trform}
A(n,3)= \frac{1}{\sqrt{n}} \sum_{Q \in {\cal Q}_{3n}/\Gamma} \chi_{n,-3}(Q) J(\alpha_Q) \, .  
\ee
If $n,d=3$ are not coprime, then according to Remark 2 following Theorem 6 of \cite{zagier} the same formula holds but with
$\chi_{n,-3}$ replaced by $0$ for imprimitive forms $Q$ which are divisible by $3$ which also divides $n$. 

\begin{exmp} (from \cite{zagier}) 
Take $n=5, d=3$. There are two $\mathrm{SL}(2,\IZ)$-equivalence classes of quadratic forms with discriminant $-15$, $Q_1=[1,1,4]$ and
$Q_2=[2,1,2]$. Their roots $\alpha_{Q_1},\alpha_{Q_2} \in \mathfrak{h}$ occur at $\alpha_{Q_1}=(1+i \sqrt{15})/2$ and $\alpha_{Q_2}=(1+i \sqrt{15})/4$ and are mapped by $J$ as
\begin{align} 
J(\alpha_1) &= (-191025-85995 \sqrt{5})/2 - 744  \\
J(\alpha_2) &= (-191025+ 85995 \sqrt{5})/2 -744 \, . \nonumber  
\end{align} 
The genus characters are $\chi_{n,-3}(Q_1)=1$
and $\chi_{n,-3}(Q_2)=-1$ so we have
\be
A(5,3)= \frac{J(\alpha_{Q_1})- J(\alpha_{Q_2})}{\sqrt{5}}=-85995 \, .
\ee
which is also the negation of the dimension of a complex conjugate pair of  irreducible representations of $\mathrm{Th}$, ${^{85995}V_{9}}$ and ${^{85995}V_{10}}$. In a manner that will be generalized in the next example, we suggestively write this in terms of characters of the identity conjugacy class of $\mathrm{Th}$ as
\be
A(5,3) = -\chi_{9}(1A)= -\chi_{10}(1A)
\ee
\end{exmp}

Alternate proofs of Zagier's results appear in \cite{bringmann} utilizing results on Maass-Poincar\'e series and their generalizations due to Niebur \cite{niebur} which have the benefit of providing
direct formulae for the Fourier coefficients of the $f_d$, as well as those of weakly holomorphic forms of weight $\lambda+\frac{1}{2}$ for
several values of $\lambda$,  in terms of Kloosterman sums. In particular, Theorem 2.1 of \cite{bringmann} with $\lambda=0$ and $m=3$ gives an expression for the
coefficients $A(n,3)$ with $n>0$ (denoted by $b_0(-3;n)$ in \cite{bringmann})
\bea \label{bkloos}
A(n,3) =  & -& 24 \delta_{\square,n} H(3) +  \pi\sqrt{2}(3/n)^{\frac{1}{4}}(1-i)  \\
& \times & \sum_{\substack{c>0 \\ c\equiv 0 ~{\rm mod}~ 4}}(1+\delta_{\mathrm{odd}}(c/4)) \frac{K_{\psi_0}(-3,n,c)}{c}I_{\frac{1}{2}}\left(\frac{4\pi\sqrt{3n}}{c}\right) \, . \nonumber
\eea
In this formula $H(d)$ is he Hurwitz-Kronecker class number with $H(3)=1/3$,
\be
\delta_{\square,m}= \begin{cases} 1 & \mbox{if m is a square} \\ 0 & \mbox{otherwise} \end{cases}
\ee
and $I_\ell(x)$ is the Bessel function of the first kind. In addition, 
for $c$ a positive multiple of 4 and $\lambda$ an integer, the weight $\frac{1}{2}$ Kloosterman sum with multiplier $\psi$ is given by
\be \label{kloos}
K_{\psi}(m,n,c) = {\sum_{d}}^*  \psi(c,d) e\left(\frac{m\overline{d}+nd}{c}\right)
\ee
where the sum runs over primitive residue classes mod $c$, $\overline{d}$ is the inverse of $d$ mod $c$,
and for any integer $k$, 
\be
\delta_{\mathrm{odd}}(k) \equiv 
\begin{cases}
 1 & k \text{ odd} \\  0  & \text{otherwise}. 
\end{cases}
\ee
Note that the standard $\theta$-multiplier $\psi_0$ is used in the definition of $A(n,3)$; its explicit form can be found in the previous section.

We now note an ambiguity that we will utilize later. In  \cite{zagier} and \cite{bringmann} the constant term in the Fourier expansion of the $f_d$ is chosen to be zero, and in particular $A(0,3)=0$. 
However adding a multiple of $\theta$ to $f_d$ does not change the modular properties or the singular terms in $f_d$. This ambiguity is reflected in the trace formulation
through the freedom to add a constant to the $J$-function. Since the genus character is trivial when $n$ is a square, adding such a constant changes only the coefficients of the
square powers of $q$ and in fact corresponds precisely to adding a multiple of $\theta$ to $f_d$. In equation (\ref{bkloos}) such a change  corresponds to changing the
coefficient of the $\delta_{\square,n}$ term. 

\subsection{Weakly holomorphic weight one-half modular forms at level N}

We now extend the results of the previous subsection to level $N$. With $\Gamma$ as defined earlier, we let $\Gamma_{\infty}$ denote the subgroup of $\Gamma$ consisting of upper-triangular matrices and define the width of $\Gamma$ at infinity to be the smallest positive integer $h$ for which $\Gamma_{\infty} = \langle T^h,-I\rangle$ where $T = \left(\begin{array}{cc} 1 & 1 \\ 0 & 1 \end{array}\right)$ and $I$ is the identity matrix. For example the width of $\Gamma_0(N)$ at infinity is easily seen to be $1$. For a given multiplier system $\psi$, we let $\alpha$ be the the real number given by 
\be
\psi\left(T^h\right) = e(\alpha).
\ee
The multiplier systems we will consider all have $\alpha = 0$. In \cite{ChengDuncanII}, Cheng and Duncan exposit a method which was pioneered by Poincar\'e and later refined by Rademacher for constructing functions symmetric under $\Gamma$.  One begins with a function of the form $q^{\mu} = e(\mu\tau)$, which is $\Gamma_{\infty}$-invariant if $h\mu+\alpha \in \IZ$, and constructs a function which is invariant under the full group $\Gamma$ by summing the images of $q^{\mu}$ under the $(\psi,w)$-action of coset representatives of $\Gamma_{\infty}$ in $\Gamma$:
\be
P^{[\mu]}_{\Gamma,\psi,w}(\tau) = \sum_{\gamma\in \Gamma_{\infty}\backslash \Gamma} q^{\mu}\vert_{\psi,w}\gamma.
\ee
In general, if $w\leq 2$ such a series does not converge locally uniformly in $\tau$, and one does not obtain a modular function holomorphic on the upper-half plane. Rademacher was thus led to attempt to regularize these Poincar\'e series by defining
\be\label{rad}
R^{[\mu]}_{\Gamma,\psi,w}(\tau) = \frac{1}{2}\delta_{\alpha,0}c_{\Gamma,\psi,w}(\mu,0) +  \lim_{K\rightarrow\infty} \sum_{\gamma\in \Gamma_{\infty}\backslash\Gamma_{K,K^2}} r^{[\mu]}_w(\gamma,\tau) q^{\mu}\vert_{\psi,w}\gamma
\ee
where
\be
\Gamma_{K,K^2} = \left\{  \left(\begin{array}{cc} a & b \\ c & d \end{array}\right)\in \Gamma\mid |c|<K, |d| < K^2   \right\}.
\ee
serves to specify the order in which the sum is taken, $\frac{1}{2}c_{\Gamma,\psi,w}(\mu,0)$ specifies a correction to the constant term in case $\alpha= 0$, and $r_{w}^{[\mu]}(\gamma,\tau)$ is a factor which regularizes the sum in the case that $w <1$, see \cite{ChengDuncanII} for exact expressions.

Indeed, this regularized expression has been proven to extend convergence to weights $w> 1$. In general convergence is poorly understood when $0\leq w \leq 1$,  but for specific cases
of $w=0$ and $w=1/2$ relevant to Monstrous moonshine and Umbral moonshine convergence has been proven in \cite{DuncanFrenkel} and \cite{ChengDuncanI} respectively. We will assume here that the particular series we deal with, all at $w=1/2$, are in fact convergent; this is supported by numerical evidence and also, as we will see in the next section, by the fact that the resulting Rademacher series are easily identified with the conjectured McKay--Thompson series of Thompson moonshine. 

One can derive Fourier expansions for these Rademacher sums in terms of Rademacher series with coefficients $c_{\Gamma,\psi,w}(\mu,\nu)$ as
\be
\label{radsum}
R^{[\mu]}_{\Gamma,\psi,w}(\tau) = q^{\mu} + \sum_{\substack{h\nu+\alpha\in \IZ \\ \nu\geq 0}}c_{\Gamma,\psi,w}(\mu,\nu)q^{\nu}
\ee
where expressions for $c_{\Gamma,\psi,w}(\mu,\nu)$ as well as a more detailed discussion of the convergence properties of Rademacher sums can be found in \cite{ChengDuncanII}.

We will now specialize to the case $\Gamma=\Gamma_0(4N)$, $w=1/2$ and $\mu=-3$ and
define a family of weakly holomorphic weight one-half forms on $\Gamma_0(4N)$ with multiplier system $\psi$ and Fourier expansion
\be 
Z^{(1,3)}_{N,\psi}(\tau) = q^{-3} + \sum_{\substack{n\geq 0 \\ n\equiv 0,1 ~{\mathrm{mod}}~ 4}} A^{(1,3)}_{N, \psi}(n) q^n \, .
\ee
These agree with the forms defined in \cite{mp} when $\psi=\psi_0$ and $N$ is odd and we have $Z^{(1,3)}_{1,\psi_0}=f_3+ 4 \theta(\tau)$. The superscripts $(1,3)$ label the function $f_3$ whose
coefficients are expressed as traces of singular moduli for $\mathrm{SL}(2,\IZ)=\Gamma_0(1)$ and anticipate possible generalizations to some of the other forms treated in \cite{zagier}. 
When $\Gamma_0(N)$ is genus zero, \cite{mp} show
that the coefficients $A^{(1,3)}_{N, \psi_0}$ can be computed either in terms of coefficients of Rademacher series as given below or in terms of traces of singular moduli.
The explicit form of the Fourier coefficients is
\be \label{zcoef}
A^{(1,3)}_{N, \psi}(0) =  4\pi\sqrt{m}(1-i) \sum_{\substack{c>0 \\ c\equiv ~{\rm mod}~ 4N}}(1+\delta_{\mathrm{odd}}(c/4))\frac{K_{\psi}(-3,0,c)}{c^{\frac{3}{2}}} 
\ee
for $n=0$ and
\be \label{nonzcoef}
A^{(1,3)}_{N, \psi}(n) =  \pi\sqrt{2}(m/n)^{\frac{1}{4}}(1-i)\sum_{\substack{c>0 \\ c\equiv 0 ~{\rm mod}~ 4N}}(1+\delta_{\mathrm{odd}}(c/4)) \frac{K_{\psi}(-3,n,c)}{c}I_{\frac{1}{2}}\left(\frac{4\pi\sqrt{3n}}{c}\right)  
\ee
for $n>0$. These expressions agree with \cite{ChengDuncanII} after projection to the Kohnen plus space and generalize the result (\ref{bkloos}) to level $N$ and generic multiplier system, but use a convention in which the constant term in the Fourier expansion is non-zero. 

We now present the generalization from \cite{mp} of (\ref{trform}) for traces of singular moduli at level $N$ when the multiplier system is the standard one, $\psi_0$. To present their result we first define the genus character of an integral binary quadratic form $Q(X,Y) = aX^2 + bXY + cY^2$ for a fundamental discriminant $D_1$ as 
\bea
\chi_{D_1}(Q) \equiv 
  \left\{\begin{array}{ll}
  0, & (a,b,c,D_1) > 1 \\
   \left(\frac{D_1}{r}\right), & (a,b,c,D_1) = 1, Q \text{ represents } r, \text{ and } (r,D_1) = 1. \end{array}\right.
\eea
Then for $D_2$ a nonzero integer with $(-1)^{\lambda}D_2 \equiv 0,1 \ \mathrm{mod} \ 4$ and $(-1)^{\lambda}D_1D_2<0$ the twisted trace is given by
\be
\mathrm{Tr}_{N,D_1}(f;D_2) \equiv \sum_{\substack{Q\in\mathcal{Q}_{|D_1D_2|}/\Gamma_0(N) \\ a\equiv 0 ~{\rm mod}~ N}} \frac{\chi_{D_1}(Q)f(\alpha_Q)}{\omega_Q}
\ee
where $\omega_Q$ is the order of the stabilizer of $Q$ in $\Gamma_0(N)/\{\pm 1\}$ \footnote{This corrects a typographical error in the description of equation (1) in \cite{mp}.}. A method for computing representatives of the space $\mathcal{Q}_d/\Gamma_0(N)$ as well as $\omega_Q$ is provided in the appendix. In terms of this trace, the coefficient of $q^{n}$ for $n>0$ and $n$ not a square in the expansion of $Z^{(1,3)}_{N, \psi_0}$ in the case that $\Gamma_0(N)$ is genus zero is given by 
\be
A^{(1,3)}_{N, \psi_0}(n) = \frac{\mathrm{Tr}_{N,-3}(T_N;n)}{\sqrt{n}} \, 
\ee
where $T_N$ is the hauptmodul for $\Gamma_0(N)$. For $n$ a square the coefficients can also be expressed in terms of traces of singular moduli, but the precise formula will in general involve adding a constant to the hauptmodul in the trace.

\begin{exmp}
Here we generalize the previous example  and illustrate the relation to traces of singular moduli at higher level by taking $N=3$ and using the hauptmodul ${T}_{3\mathrm{B}}(\tau) = \eta(\tau)^{12}/\eta(3\tau)^{12} + 12$. To compute the coefficient of $q^5$ in $Z^{(1,3)}_{3, \psi_0}$, we sum over representatives of the space $\mathcal{Q}_{15}/\Gamma_0(3)$. Again, using standard reduction theory we find that every binary quadratic form in $\mathcal{Q}_{15}$ is equivalent to $[1,1,4]$ or $[2,1,2]$ under the action of $\mathrm{SL}(2,\IZ)$, and it is easily deduced that a complete set of representatives of the quotient $\mathcal{Q}_{15}/\Gamma_0(3)$ can be found in the set
\be 
\mathcal{R} = \Big\{ Q|_{\gamma} : Q = [1,1,4], [2,1,2],\ \gamma \in \Gamma_0(3) \Big\}.
\ee
The elements of $\mathcal{R}$ with $a \equiv 0 \ \mathrm{mod} \ 3$ are $Q_1 = [3,-3,2]$ and $Q_2 = [6,3,1]$ and these are inequivalent under the action of $\Gamma_0(3)$. These are trivially stabilized in $\Gamma_0(3)$
and have 
\be
\chi_{-3}(Q_1) = -1, \ \ \ \chi_{-3}(Q_2) = 1  
\ee
so that
\be
A^{(1,3)}_{3, \psi_0}(5) = \frac{{ T}_{3\mathrm{B}}\left(\alpha_{Q_2} \right)-T_{3\mathrm{B}}\left(\alpha_{Q_1}\right)}{\sqrt{5}} = 27 \, .
\ee
We note that this is the negative of the character of the $3$B conjugacy class in $\mathrm{Th}$ for both irreducible representations of dimension $85995$:
\be
A^{(1,3)}_{3, \psi_0}(5)  = -\chi_{9}(3B) = - \chi_{10}(3B). 
\ee
A more detailed method for the computation of these coefficients is provided in Appendix \ref{appa}. 
\end{exmp}
The previous two examples highlight a possible connection between the coefficients of $Z^{(1,3)}_{N, \psi_0}$ and twisted versions ${\cal F}_3$ obtained by replacing its coefficients, thought of as characters of the identity class, with characters of other conjugacy classes. Indeed, this observation allows us to obtain the McKay--Thompson series of $\mathcal{F}_3$ by considering more general multiplier systems as detailed in the next section.

\section{McKay--Thompson series for the Thompson group}

As in Monstrous moonshine and Umbral moonshine, after associating supermodules to the Fourier-coefficients in the $q$-expansion of $\mathcal{F}_3$, we can twist  $\mathcal{F}_{3}$ by Th conjugacy classes $[g]$ to obtain the McKay--Thompson series
\be{\cal F}_{3,[g]}(\tau)= \sum_{\substack{m=-3 \\ m \equiv 0,1 ~{\rm mod}~ 4}}^{\infty}\mathrm{str}_{W_m}(g) q^m .
\ee
Our main goal in this section is to exhibit a  weight one-half weakly holomorphic modular form for each $[g]$ and decompositions for the components of each $W_m$ into irreducible representations of
Th with positive integer multiplicity such that the above trace function recovers the corresponding modular form. 

As a simple example consider the $3$B McKay--Thompson series with leading terms 
\be
\mathcal{F}_{3,3\mathrm{B}}(\tau) = 2q^{-3} + 5+ 54q^{4}+ 54q^{5}-108q^{8} + 16q^{9}+ 12q^{12}+ \cdots
\ee
where we have used the decompositions of coefficients of $q^m$ for $m\leq 12$ appearing in Table \ref{thomqk}. 
The Rademacher series at level $N=3$ and multiplier system $ \psi_0$ is given by 
\be
Z^{(1,3)}_{3,\psi_0}(\tau) = q^{-3} -1/2 -5q + 22q^4 + 27q^5- 54q^8 + 3q^9 + 6q^{12} + \cdots \nonumber
\ee
and since the group $\Gamma_0(3)$ is genus $0$ with hauptmodul $T_3$  we can compute traces of singular moduli and find agreement amongst the coefficients of $q^m$ with $m$ not a square. 
 As mentioned earlier, the coefficients of square powers of $q$ depend on the
choice of constant term in the hauptmodul in the trace formulation and this is reflected in the Rademacher series by the possibility to add terms involving theta functions without
changing the modular properties. We can use this to make the identification
\be
\mathcal{F}_{3,3\mathrm{B}}(\tau) = 2Z^{(1,3)}_{3, \psi_0}(\tau)+ 6 \theta(\tau)
\ee
which suggests that $\mathcal{F}_{3,3\mathrm{B}}$ is modular for  $\Gamma_0(12)$. Computation to higher order in the $q$-expansion confirms this identification. 

Although \cite{mp} do not prove that the trace formulation agrees with the Rademacher series coefficients for even values of $N$, we find nonetheless that 
similar computations lead to agreement between the trace formulation and Rademacher coefficients for $Z^{(1,3)}_{4, \psi_0}(\tau)$ and $Z^{(1,3)}_{8, \psi_0}(\tau)$ and that these can be matched onto the $4$A and $4$B twists of $\mathcal{F}_3$ with Fourier expansions
\begin{align}
\mathcal{F}_{3,4\mathrm{A}}(\tau) &= 2q^{-3} + 8 + 16q^4 + 42q^{5}-84q^{13}+ \cdots \\
\mathcal{F}_{3,4\mathrm{B}}(\tau) &= 2q^{-3} -22q^{5}+108q^{13}+ \cdots
\end{align}
as
\be
\mathcal{F}_{3,4\mathrm{A}}(\tau) = 2\big(Z_4(\tau)+ 4 \theta(4\tau)\big)
\ee
\be
\mathcal{F}_{3,4\mathrm{B}}(\tau) = -2\big(Z_4(\tau)-2Z_8(\tau)\big).
\ee
Note that because 4 is a square, we are free to add a multiple of $\theta(4\tau)$ to a weakly holomorphic weight one-half modular form over $\Gamma_0(16)$ or $\Gamma_0(32)$ while retaining the Kohnen plus condition on its Fourier coefficients and without altering its modular properties. 

Continuing in this way, we find relations of the form
\be
{\cal F}_{3,[g]}(\tau)= 2Z^{(1,3)}_{o(g),\psi_0}(\tau) + \sum_{\substack{m > 0 \\ m^2|o(g)}} \kappa_{m,g} \theta\left(m^2 \tau\right)
\ee
for the conjugacy classes $[g]$ 
of Th labelled as 2A, 3B, 4A, 5A, 6C, 7A, 9A, 9B, 10A, 12C, 13A, and 18A in Tables \ref{thchartabone}--\ref{thchartabfour}. The constants
$\kappa_{m,g}$ are as specified in Table \ref{theconst}. In addition we have
\bea
{\cal F}_{3,4\mathrm{B}}(\tau) &=& -2\big(Z_4(\tau) - 2 Z_8(\tau)\big)\\
{\cal F}_{3,8\mathrm{A}}(\tau) &=& -2\big(Z_8(\tau) - 2Z_{16}(\tau)\big) \, .
\eea
Furthermore, we find for these cases, all of which involve $o(g)$ such that $\Gamma_0(o(g))$ is genus zero and whose hauptmoduls are tabulated in Table \ref{hauptmoduln}, that the trace formulations and Rademacher coefficients agree for coefficients of non-square powers of $q$. 

We mentioned earlier that the results of \cite{mp} that we have used do not strictly apply for
$N$ even. However the fact that very similar results were shown to extend to $N=2$ in \cite{benlar}, the fact that the results for even $N$
when $\Gamma_0(N)$ is genus zero coincide with the McKay--Thompson series above, and the fact that the results of \cite{mp} for these values of $N$ coincide with
numerical results obtained from the computation of coefficients of  Rademacher series convince
us that it is very likely that the results of \cite{mp} can be extended at least to even $N$ for which $\Gamma_0(N)$ is genus zero. 

\begin{table}\begin{center}
\setlength{\tabcolsep}{.8cm}
\begin{tabular}{llr}
\toprule
Monster class & $\Gamma_0(N)$ & $T_N$\\ \toprule
2B & $\Gamma_0(2)$ & $\eta(\tau)^{24}/\eta(2\tau)^{24} + 24$ \\
3B & $\Gamma_0(3)$ & $\eta(\tau)^{12}/\eta(3\tau)^{12} + 12$ \\
4C & $\Gamma_0(4)$ & $\eta(\tau)^8/\eta(4\tau)^8 +  8$ \\
5B & $\Gamma_0(5)$ & $\eta(\tau)^6/\eta(5\tau)^6 + 6$ \\
6E & $\Gamma_0(6)$ & $\eta(\tau)^5\eta(3\tau)/\eta(2\tau)\eta(6\tau)^5+5$ \\
7B & $\Gamma_0(7)$ & $\eta(\tau)^4/\eta(7\tau)^4+4$ \\
8E & $\Gamma_0(8)$ & $\eta(\tau)^4\eta(4\tau)^2/\eta(2\tau)^2\eta(8\tau)^4+4$ \\
9B & $\Gamma_0(9)$ & $\eta(\tau)^3/\eta(9\tau)^3+3$ \\
10E & $\Gamma_0(10)$ & $\eta(\tau)^3\eta(5\tau)/\eta(2\tau)\eta(10\tau)^3+3$ \\
12I & $\Gamma_0(12)$ & $\eta(\tau)^3\eta(4\tau)\eta(6\tau)^2/\eta(2\tau)^2\eta(3\tau)\eta(12\tau)^3+3$\\
13B & $\Gamma_0(13)$ & $ \eta(\tau)^2/\eta(13\tau)^2+2$ \\
16B & $\Gamma_0(16)$ & $\eta(\tau)^2\eta(8\tau)/\eta(2\tau)\eta(16\tau)^2+2$ \\
18D & $\Gamma_0(18)$ & $\eta(\tau)^2\eta(6\tau)\eta(9\tau)/\eta(2\tau)\eta(3\tau)\eta(18\tau)^2+2$\\
(25Z) & $\Gamma_0(25)$ & $\eta(\tau)/\eta(25\tau) + 1$\\\bottomrule
\end{tabular}
\caption{\label{hauptmoduln} The McKay--Thompson series for the monster which are hauptmoduls for genus 0 congruence subgroups, $\Gamma_0(N)$. The hauptmoduls are normalized to have
leading singularity $q^{-1}$ at the cusp at infinity and vanishing constant term.}\label{table:hauptmoduls}
\end{center}
\end{table}

This construction involving traces of singular moduli for genus zero $\Gamma_0(N)$ and Rademacher series with multiplier $\psi_0$  is not sufficient to generate candidate forms for all the McKay--Thompson series of Thompson moonshine. However the expressions for the coefficients in equations (\ref{zcoef}), (\ref{nonzcoef})
can be modified to  obtain weight one-half weakly holomorphic modular forms at level $4 N h$ by simply changing the multiplier system $\psi$ used in the  Kloosterman sums $K_\psi(-3,n,c)$. We will work with multiplier systems of the form
\be \label{newmult}
\psi_{N,v,h}(\gamma)=\psi_0(\gamma)e\left(-  v\frac{c d}{Nh}\right)
\ee
with $v, h \in \IZ$ and $h$ dividing $4 \times 24$, see \cite{ChengDuncanII} for a general discussion including the fact that the last factor in (\ref{newmult}) is actually a character for
$\Gamma_0(4N)$.

We thus assign to each element $g$ of Th a pair $(v_g,h_g)$ depending only on the class $[g]$ with the requirement that $h_g$ divides $2o(g)$ and find that we can make
the identification
\be
{\cal F}_{3,[g]}(\tau)= 2 Z^{(1,3)}_{o(g),\psi}(\tau) + \sum_{\substack{m>0 \\  m^2| o(g) h_g}} \kappa_{m,g}\theta\left(m^2 \tau\right), \ \ \ \psi = \psi_{o(g),v_g,h_g}
\ee
where the quantity on the right hand side is, modulo issues of convergence, a weakly holomorphic weight one-half modular form on $\Gamma_0(4 o(g) h_g)$.
The pairs $(v_g,h_g)$ specifying the multiplier systems as well as the $\kappa_{m,g}$ may be found in Table \ref{theconst}.

Given this identification of the McKay--Thompson series we can solve for the multiplicities $b^{(a)}_{m,i}$ that specify the multiplicity of
the representation $V_i$ in the decomposition of $W_m^{(a)}$. These are tabulated in Tables \ref{decomp1} and \ref{decomp2}
for $m \le 32$. We have in fact verified that the $b^{(a)}_{m,i}$ are positive integers for $m \le 52$. 

We summarize the discussion here and in the previous sections with the following.
\begin{conj}
We conjecture that there exists a naturally defined $\IZ$--graded supermodule for the Thompson sporadic group
\be
W= \bigoplus_{\substack{m \ge -3 \\ m \equiv 0,1~ {\rm mod} ~4}}^\infty W_m \, ,
\ee
with $W_m$ for $m \ge 0$ and $m=0 ~{\mathrm{mod}}~ 4$ having vanishing odd part,
$W_m$ for $m \ge 0$ and  $m=1 ~{\mathrm{mod}}~4$ having vanishing even part, and $W_{-3}$ having vanishing odd part,
such that the graded dimensions of $W_m$ are related to the weight one-half weakly holomorphic modular form ${\cal F}_3$
by
\be
{\cal F}_3(\tau)=  \sum_{\substack{m=-3 \\ m \equiv 0,1 ~{\rm mod}~ 4}}^{\infty} \mathrm{str}_{W_m}(1) q^m \, ,
\ee
and such that the weight one-half weakly holomorphic modular forms $\mathcal{F}_{3,[g]}$
described above and in Table \ref{mtseries} are related to $W$ via graded supertrace functions via
\be
{\cal F}_{3,[g]}(\tau)=  \sum_{\substack{m=-3 \\ m \equiv 0,1 ~{\rm mod}~ 4}}^{\infty} \mathrm{str}_{W_m}(g) q^m \, .
\ee
\end{conj}

\section{A discriminant property for Thompson moonshine}

We now describe a discriminant property of the moonshine connection between the Thompson group and  ${\cal F}_3$.
We start with a quick review of the discriminant property of Umbral moonshine, concentrating for simplicity on the cases where
the Niemeier lattice $X$ is of pure A-type, 
$X=A_{\ell-1}^{24/(\ell-1)} $,
with $\ell-1$ a divisor of $12$ as discussed in  \cite{um}.

The main ingredients for the A-type examples of Umbral moonshine are $\ell-1$ component vector-valued mock modular forms $H^{(\ell)}$ and a finite
group 
\be
G^{(\ell)}= \mathrm{Aut}\big(L^{(\ell)}/W^{(\ell)}\big)
\ee
with $L^{(\ell)}$ the Niemeier lattice with root system $A_{\ell-1}^{24/(\ell-1)}$ and $W^{(\ell)}$
the Weyl group of $L^{(\ell)}$ generated by reflections in the roots. The r-th component $H^{(\ell)}_r$ of $H^{(\ell)}$ has a Fourier expansion of the
form 
\be
H^{(\ell)} _r(\tau)= \sum_{k \in \IZ} c_r^{(\ell)}(k-r^2/4 \ell) q^{k-r^2/4 \ell}
\ee
with $c_r^{(\ell)}(k-r^2/4 \ell)$ the dimension of a $G^{(\ell)}$ module $K^{(\ell)}_{r,k-r^2/4 \ell}$. One says that an integer $D$ is a discriminant of $H^{(\ell)}$ if
there exists a term $q^{-D/4 \ell}$ with non-vanishing Fourier coefficient in at least one of the components of $H^{(\ell)}$.

The main results relating the number fields over which the irreducible representations of $G^{(\ell)}$ are defined and the discriminants of $H^{(\ell)}$ 
are stated in Propositions 5.7 and 5.10 of \cite{um} which we state here for completeness.

\begin{proposition}\label{umdiscone}
{\rm (Proposition 5.7 of \cite{um})}
If $n>1$ is an integer satisfying 
\begin{enumerate}
\item{there exists an element of $G^{(\ell)}$ of order $n$}, and
\item{there exists an integer $\lambda$ that is co-prime to $n$ such that $D = -n \lambda^2$ is a discriminant of $H^{(\ell)}$,}
\end{enumerate}
then there exists at least one pair of irreducible representations $\varrho$ and $\varrho^*$ of $G^{(\ell)}$ and at least one element $g \in G^{(\ell)}$ such that $\tr_{\varrho}(g)$ is not rational but
\be\label{n_type}
{\tr}_{\varrho} (g), {\tr}_{ \varrho^*} (g) \in \IQ(\sqrt{-n})
\ee
and $n$ divides $o(g)$.
\end{proposition}

A representation $\rho$ of $G^{(\ell)}$ is said to be of {\it type $n$} if $n$ is an integer satisfying the two conditions of Proposition \ref{umdiscone} and the character
values of $\rho$ generate $\IQ[\sqrt{-n}]$.  A connection between discriminants of $H^{(\ell)}$ and representations of type $n$ is then described by
\begin{proposition}\label{umdisctwo}
{\rm (Proposition 5.10 of \cite{um})}
Let $n$ be one of the integers satisfying the two conditions of Proposittion \ref{umdiscone}  and let $\lambda_{n}$ be the smallest positive integer such that $D = -n \lambda_{n}^2$ is a discriminant of $H^{(\ell)}$. Then $K^{(\ell)}_{r,-D/4 \ell }= \varrho_{n} \oplus \varrho_{n}^*$ where $\varrho_{n}$ and $ \varrho_{n}^*$ are dual irreducible representations of type $n$. Conversely, if $\varrho$ is an irreducible representation of type $n$ and $-D$ is the smallest positive integer such that $K^{(\ell)}_{r,-D/4 \ell}$ has $\varrho$ as an irreducible constituent then there exists an integer $\lambda$ such that $D = - n \\lambda^2$. 
\end{proposition} 

The latter was extended to the following two conjectures in \cite{um}
\begin{conj}\label{conj:conj:disc:dualpair}
{\rm (Conjecture 5.11 of \cite{um} )}
If $D$ is a discriminant of $H^{(\ell)}$ which satisfies $D = -n \lambda^2$  for some integer $\lambda$ then the representation $K^{(\ell)}_{r,-D/4 \ell}$ has at least one dual pair of irreducible representations of type $n$ arising as irreducible constituents. 
\end{conj}
\begin{conj} {\rm Conjecture 5.12 of \cite{um} )}
For $\ell \in \Lambda=\{2,3,4,5,7,13\}$ the representation $K^{(\ell)}_{r,-D/4 \ell} $is a doublet if and only if $D \ne n \lambda^2$ for any integer $\lambda$ for any n satisfying the conditions of Proposition 5.7.
\end{conj}
In the above a $G$-module $V$ is a {\it doublet} if it is isomorphic to the direct sum of two copies of a single representation of $G$.

For $\ell=2$ these conjectures are proved in \cite{Creutzig}.  Further details of the discriminant structure and the extension to general Niemeier root systems $X$ can be found
in \cite{um,mum}.

We now describe an analogous structure in Thompson moonshine. 
We first note that it follows from the results of \cite{zagier}, summarized in section 2,  that the coefficient $c(m)$ of $q^m$ in ${\cal F}_3$ can be computed in terms of traces of
singular moduli which are determined by the root of a quadratic form of discriminant $-3m$ (i.e. a Heegner point of discriminant $-3m$). We will say that $-3m$ is a discriminant of ${\cal F}_3$ if there exists a term $q^m$ in the Fourier
expansion of ${\cal F}_3$ with non-zero coefficient. 
Any discriminant can be written uniquely in the form $D_0 \lambda^2$ with $\lambda$
a positive integer and $D_0$ a fundamental discriminant. Thus for each $m >0$ with $m= 0,1 ~{\rm mod}~ 4$ we have a unique decomposition $-3m= D_0(m) \lambda^2$ with $D_0(m)$ a negative 
fundamental discriminant. 

We now state a proposition relating the discriminants of ${\cal F}_3$ to properties of the characters of Th.
\begin{proposition} \label{discone}
If $D_0(m)$ is a negative fundamental discriminant satisfying
\begin{enumerate}
\item there exists an element of {\rm Th} of order $|D_0(m)|$, and
\item there exists a positive integer $\lambda$ such that $-3m= D_0(m) \lambda^2$ is a discriminant of ${\cal F}_3$ and $(\lambda,3)=1$,
\end{enumerate}
then there exists at least one pair of irreducible representations $V$ and $\overline{V}$ of {\rm Th} and at least one element $g \in Th$ such that $\mathrm{tr}_V(g)$ is not rational
but
\be
\mathrm{tr}_{V}(g), \mathrm{tr}_{\overline{V}}(g) \in \IQ\Big[\sqrt{D_0(m)}\Big]
\ee
and $|D_0(m)|$ divides $o(g)$.
\end{proposition}

The negative fundamental discriminants obeying the two conditions of Proposition \ref{discone} are
$-3$, $-15$,$-24$, and $-39$ and inspection of the character table of Th shows that Proposition \ref{discone} is true since the pairs of irreducible representations
${}^{27000}V_4,{}^{27000}V_5$, ${}^{85995}V_9, {}^{85995}V_{10}$, ${}^{1707264}V_{17}, {}^{1707264}V_{18}$ and ${}^{779247}V_{14},{}^{779247}V_{15}$
have characters in $\IQ\big[\sqrt{D_0}\big]$ for elements $g$ of Th with $|D_0|$ dividing $o(g)$ for $D_0=-3, -15,-24,-39$ respectively.  

We will say that
a representation $V$ of Th is of type $D_0$ if $D_0$ is a negative fundamental discriminant satisfying the two conditions of Proposition (\ref{discone}) and the
character values of $V$ generate the ring of algebraic integers in $\IQ\big[\sqrt{D_0}\big]$ over $\IZ$. The latter condition excludes representations defined over
$\IQ\big[ \sqrt{-3} \big]$ and leaves the list of irreducible representations of type $D_0$ is given in Table \ref{typerep}.
\begin{table} 
\begin{center}
\begin{tabular}{cc}\toprule
$D_0$ & $(V, \overline{V})$  \\\midrule  
$-15$ & $(V_9,V_{10}),(V_{35},V_{36}) $\\
$-24$ & $(V_{17},V_{18})$  \\
$-39$ & $(V_{14},V_{15})$  \\
\bottomrule
\end{tabular}
\end{center}\caption{\label{typerep} {Irreducible representations of Th of type $D_0$.}}
\end{table}

There is also a connection between the discriminants of ${\cal F}_3$ and the fields over which the representations $W_m$ attached to the coefficient of $q^m$
in ${\cal F}_3$ are defined which we state as
\begin{proposition} \label{disctwo}
Let $D_0(m)$ be one of the fundamental discriminants satisfying the two conditions of Proposition \ref{discone} and let $\lambda_m$ be the smallest positive
integer such that $-3m=D_0(m) \lambda_m^2$ is a discriminant of ${\cal F}_3$. Then $W_m= V \oplus \overline{V}$ where $V$ and $\overline{V}$ are dual irreducible representations
defined over  $\IQ\big[ \sqrt{D_0(m)} \big]$. 
\end{proposition}

\begin{rmk}
Since the Schur index of all irreducible representations of {\rm Th} is one \cite{feit}, it follows that {\rm Th} representations of type $D_0$ can be realized over $\IQ\big[\sqrt{D_0}\big]$.
\end{rmk}

Conjecture 5.11 of \cite{um} also predicts properties of the representations attached to discriminants of the mock modular forms $H^X$ to arbitrary powers in
their Fourier expansion.   The natural generalization of Conjecture 5.11 to our situation is
\begin{conj}
Whenever $-3m$ is a discriminant of ${\cal F}_3$
with $-3m = D_0(m) \lambda^2$ with $|D_0(m)|$ an order of Th and $(\lambda,3)=1$ then the representation $W_m$ has at least one pair of irreducible representation of
type $D_0$ arising as irreducible constituents. 
\end{conj}
Evidence for this conjecture can be seen from an inspection of the multiplicities of representations appearing in Tables \ref{decomp1} and \ref{decomp2}.

Finally we have a conjecture that is similar in nature to conjecture 5.12 of \cite{um} regarding the presence of decompositions with odd multiplicities of representations of
type $D_0(M)$. 
\begin{conj}
If $-3m$ is a discriminant of ${\cal F}_3$
with $-3m = D_0(m) \lambda^2$ with  $(\lambda,3)=1$  and $D_0(m) \in \{-15,-24,-39 \}$ then the representation $W_m$ contains a pair $V \oplus \bar V$ of irreducible representations of
type $D_0$  with odd multiplicity.
\end{conj}

\begin{rmk}
According to Corollary 1.2 of \cite{ono}, the coefficients of $q^{-D/4 \ell}$ appearing in the discriminant  property of Umbral Moonshine for $X$ of pure $A$-type are also related to traces of singular moduli 
evaluated at Heegner points of discriminant $D$.
\end{rmk}

\paragraph{Acknowledgments}

This work was inspired in part by a question asked by T. Piezas III on {\tt mathoverflow.net} in which he noted a connection between
the Thompson group and the weight one-half modular form $f_3 + 4 \theta(\tau)$.  We thank him for asking the question, S. Carnahan for useful comments in response to his question, and {\tt mathoverlow.net}  for providing a forum for such questions. We owe particular thanks to J. Duncan for remarks that led to substantial improvements in the manuscript.  JH thanks the Institut 
\mbox{d'\'Etudes} Scientifiques de Carg\`ese for hospitality during an early stage of this work. BR thanks Izaak Meckler for his insight on implementing computations in software and more generally for his willingness to discuss various aspects of moonshine. This work was supported in part by NSF Grant 1214409. 

\vfil \eject

\appendix

\section{Computing traces of singular moduli} \label{appa}

We present methods for computing representatives for the quotient space $\mathcal{Q}_d/\Gamma_0(N)$ as well as a generalized reduction theory of integral binary quadratic forms, which is similar to the one found in \cite{gcb}. Our algorithms rely on knowledge of $\mathcal{Q}_d/\mathrm{SL}(2,\IZ)$ from reduction theory and the structure of $\Gamma_0(N)\backslash\mathrm{SL}(2,\IZ)$.

The basic idea is this. Let $X$ be a well-ordered set on which a group $G$ acts from the left with finitely many orbits and $H$ be a subgroup of $G$ with finite index. We assume that the stabilizer of any element in $X$ is finite. Then once one fixes a set of orbit and right coset representatives, $(X/G)^{rep}$ and $(H\backslash G)^{rep}$, the finite set
\be
\mathcal{R} = \{ g\cdot x \mid x\in (X/G)^{rep}, \ g \in (H\backslash G )^{rep}  \}
\ee
contains a complete set of representatives for the space $X/H$. However, the existence of non-trivial stabilizers associated to the group action of $G$ implies that there may be two elements in $\mathcal{R}$ which are equivalent under $H$. To this end, we will prescribe a method for testing whether or not two elements are equivalent under the action of $H$ and filter elements out of $\mathcal{R}$ until we have a true set of representatives for the quotient $X/G$. 

Eventually, we will apply this to $X = \mathcal{Q}_d$, $G=\mathrm{SL}(2,\IZ)$, and $H=\Gamma_0(N)$, and to this end, we state some facts from computational group theory. These results also hold for right actions and left coset representatives.

\subsection{Reducing to fundamental domain of subgroup}

Fix a fundamental domain $(X/G)^{rep}$ and a set of right coset representatives $(H\backslash G)^{rep}$. Define the reduction function
\begin{align}
R_G:X&\rightarrow X\times G \\
x&\mapsto (x_G,g_x)
\end{align}
which computes the \emph{unique} representative of $x$ in $(X/G)^{rep}$ as well as some group element which connects them, $g_x\cdot x = x_G$. Note in general that $g_x$ is \emph{not} unique. If $x_G$ is stabilized by $\{s_1,\dots,s_m\}$, then the reduction function could have returned $s_i g_x$ for $i = 1,\dots, m$ as the group element connecting $x$ to $x_G$. In fact, one can see that these are the only group elements connecting $x$ to its reduction $x_G$. We thus give an algorithm for computing $R_H$ in terms of $R_G$ which copes with the existence of non-trivial stabilizers.

Let $(x_G,g_x) = R_G(x)$ and Stab$_G(x_G) = \{s_1,\dots,s_m\}$. Then for each $k=1,\dots, m$, it is easy to see that there is a unique coset representative $g_{k}\in (H\backslash G)^{rep}$ for which $g_ks_kg_x \in H$. We can use the ordering on $X$ to define $h = g_ks_kg_x$ using the $k$ which minimizes $(g_ks_kg_x)\cdot x$. Then we let
\be
R_H(x) = (h\cdot x,h).
\ee
This definition of $R_H$ ensures that every pair of elements $x$ and $y$ in $X$ which are equivalent under the action of $H$ are mapped onto the \emph{same} reduced element, by perscribing a canonical way to choose amongst the $m$ different group elements $s_kg_x$ which connect $x$ to $x_G$. 
\begin{exmp}
Let $X = \mathcal{Q}_3$ with lexicographic order, defined as $Q_1 = [a_1,b_1,c_1] < Q_2 = [a_2,b_2,c_2]$ if $a_1<a_2$ or $a_1=a_2$ and $b_1<b_2$ or $a_1=a_2,\ b_1=b_2$ and $c_1<c_2$. Further let $G$ the full modular group $\mathrm{SL}(2,\IZ)$, and $H$ the congruence subgroup $\Gamma_0(2)$. We would like to compute $R_{\Gamma_0(2)}(Q)$ where $Q = 3X^2 + 3XY + Y^2$. 

From reduction theory, we get that 
\be
R_{\mathrm{SL}(2,\IZ)}(Q) = (Q_{\mathrm{SL}(2,\IZ)},M)\equiv\left(X^2 + XY + Y^2,\left(\begin{array}{rr}0 & -1 \\ 1 & 2 \end{array}\right)\right).
\ee 
We can choose left coset representatives 
\be
(\Gamma_0(2)\backslash{\mathrm{SL}(2,\IZ)})^{rep} = \left\{M_1,M_2,M_3\right\}= \left\{ \left(\begin{array}{cc} 1 & 0 \\ 0 & 1 \end{array}\right),  \left(\begin{array}{cr} 0 & -1 \\ 1 & 0 \end{array}\right),  \left(\begin{array}{cc} 1 & 0 \\ 1 & 1 \end{array}\right)  \right\}
\ee
and note that 
\be
\mathrm{Stab}_{{\mathrm{PSL}(2,\IZ)}}(Q_{{\mathrm{SL}(2,\IZ)}}) = \{S_1,S_2,S_3\} =\left\{ \left(\begin{array}{rr} 1 & 0 \\ 0 & 1 \end{array}\right),  \left(\begin{array}{rr} -1 & -1 \\ 1 & 0 \end{array}\right),  \left(\begin{array}{rr} 0 & 1 \\ -1 & -1 \end{array}\right) \right\}.
\ee
Then, the matrices $MS_1M_2$, $MS_2M_3$, and $MS_3M_1$ all lie in $\Gamma_0(2)$. Acting with these matrices from the right on $Q$ gives the binary quadratic forms $X^2-XY+Y^2$, $3X^2 + 3XY + Y^2$, and $X^2+XY+Y^2$ respectively. The least of these lexicographically is the first, and so we have that 
\be
R_{\Gamma_0(2)}(3X^2+3XY+Y^2) = \left(X^2-XY+Y^2,\left(\begin{array}{rr}1 & 0 \\ 2 & 1 \end{array}\right)\right).
\ee
\end{exmp}
\subsection{Representatives for quotient by a subgroup}
It is easy to see that once we have a true ``reduction algorithm" we can easily test whether or not two elements are equivalent under the action of the group.

\begin{claim}
For two elements $x$ and $y$ of $X$, let $(x_h,g_x) = R_H(x)$ and $(y_h,g_y) = R_H(y)$. These are equivalent under the action of $H$ if and only if $x_h=y_h$.
\end{claim}

\begin{rmk}
This claim also allows us to implicitly define a fundamental domain of $X$ under the group action of $H$ as $(\pi_1\circ R_H)(X)$ where $\pi_1:X\times G\rightarrow X$ denotes the projection onto the first factor. 
\end{rmk}

Now that we are able to test equivalence, we can carry out the procedure suggested at the beginning of this section.

\begin{exmp}
We would like to compute the space $\mathcal{Q}_4/\Gamma_0(3)$. From reduction theory we know that every form with discriminant $-4$ is equivalent under the action of the full modular group to $X^2 + Y^2$. Acting on this form with the coset representatives of $\Gamma_0(3)\backslash{\mathrm{SL}(2,\IZ)}$, we get that a full set of representatives of $\mathcal{Q}_4/\Gamma_0(3)$ must be contained in 
\be
\{ X^2 + Y^2, 2X^2-2XY + Y^2, 5X^2 + 4XY + Y^2  \}.
\ee
However, we find that 
\be
R_{\Gamma_0(3)}(5X^2+4XY + Y^2) = 5X^2+4XY +Y^2 = R_{\Gamma_0(3)}(2X^2-2XY+Y^2)
\ee
while $X^2+Y^2$ and $2X^2-2XY+Y^2$ are inequivalent under $\Gamma_0(3)$, and so we have that 
\be
\big(\mathcal{Q}_4/\Gamma_0(3)\big)^{rep} = \{X^2+Y^2,2X^2-2XY+Y^2\}.
\ee
\end{exmp}
\vfil \eject

\section{Tables}

In this appendix we provide tables giving the multiplier systems and coefficients of theta functions used in specifying the McKay--Thompson series, the Fourier expansions
of all the McKay--Thompson series ${\cal F}_{3,g}$ up to order $q^{32}$,  the character table of the Thompson group and a table providing the decomposition into irreducible
representations of the modules $W_m^{(0)}$ and $W_m^{(1)}$. The character table was computed using the GAP4 computer algebra
package \cite{GAP4}. 

\begin{sidewaystable}
\begin{center}
\begin{small}
\smallskip
\begin{tabular}{ccccccccccc}
\toprule
$[g]$	&1A	&2A	&3A	&3B	&3C	&4A	&4B	&5A	&6A	&6B		\\
	\midrule
$\nu,h$ &$0,1$  &$0,1$ &$2,3$ &$0,1$ &$1,3$ &$0,1$ &$1,8$ &$0,1$ &$1,6$ & $1,3$   \\
	\midrule
$\kappa_{m,g}$ &$240_1$ &$0$ &-$6_1+18_9 $&$6_1$ &$ 0$&$8_4$ &$ 0$&$ 0$&$0 $&$ 0 $\\
\midrule 
\noalign{\vskip 10mm}
\midrule
$[g]$	&6C &7A & 8A & 8B & 9A & 9B & 9C & 10A & 12AB & 12C  	\\
\midrule
$v,h$ &$0,1$  &$0,1$ &$1,8$ &$3,16$ &$0,1$ &$0,1$ &$2,3$ &$0,1$ &$5,12$ & $0,1$    \\
	\midrule
$c_k$ &$0 $ &$2_1$ &$0 $ &$0 $ &$ 6_9$ &-$3_9$ &$0 $ &$0 $& -$1_4$+$3_{36}$ & -$1_4$   \\
\midrule
\noalign{\vskip 10mm}
\midrule
$[g]$	&  12D & 13A & 14A & 15AB & 18A & 18B & 19A &20A & 21A & 24AB  	\\
\midrule
$v,h$ &$5,24$  &$0,1$  &$0,1$ &$2,3$ &$0,1$ &$1,3$ &$0,1$ &$1,8$ &$2,3$ &$5,24$       \\
	\midrule
$\kappa_{m,g}$ &$ 0$ &$(\frac{1}{3})_1$ &$0 $ &$0 $ &$0 $ &$0 $ &$(\frac{3}{5})_1$ &$0 $ &$(\frac{1}{8})_1$-$(\frac{3}{8})_9 $ & $0 $   \\
\midrule
\noalign{\vskip 10mm}
\midrule
$[g]$ & 24CD & 27A & 27BC & 28A & 30AB & 31AB & 36A & 36BC & 39AB &  \\
\midrule
$v,h$ &$11,48$  &$2,3$ &$2,3$ &$0,1$ &$1,3$ &$0,1$ &$0,1$ &$0,1$ &$2,3$ &      \\
	\midrule
$\kappa_{m,g}$ &$ 0$ &-$1_9$+$3_{81}$ &$(\frac{1}{2})_9$-$(\frac{3}{2})_{81}$ &$1_4$ &$0$ &-$(\frac{1}{4})_1 $&$2_4$-$3_{36}$ & -$1_4$ &-$(\frac{3}{7})_1$+$(\frac{9}{7})_9$  &   \\
\\\bottomrule
\end{tabular}
\caption{Multipliers and theta function corrections for the McKay Thompson series. The notation $\sum_m (\kappa_{m,g})_m$ indicates the addition of a term $\sum_m \kappa_{m,g} \theta(m^2 \tau)$
to the modular form associated with the conjugacy class $[g]$. The
parentheses are omitted in case $\kappa_{m,g}$ is an integer and a zero entry indicates that all $\kappa_{m,g}=0$.} \label{theconst}
\smallskip
\end{small}
\end{center}
\end{sidewaystable}

\begin{sidewaystable}
\begin{small}
\centering
\caption{Coefficient of $q^n$ in the McKay--Thompson series ${\cal F}_{3,[g]}(\tau)$ of Thompson moonshine, part one.}\label{mtseries}\smallskip
\begin{tabular}{c@{ }@{\;}r@{ }r@{ }r@{ }r@{ }r@{ }r@{ }r@{ }r@{ }r@{ }r@{ }r@{ }r@{ }r@{ }r@{ }r@{ }r@{ }r@{ }r@{ }r@{ }r@{ }r@{ }r@{ }r@{ }r@{ }r@{ }r}\toprule
$n ~\backslash ~[g]$	&1A	&2A	&3A	&3B	&3C	&4A	&4B	&5A	&6A	&6B	&6C	&7A	&8A	&8B	&9A	&9B	&9C	&10A	&12A  & 12B	&12C   & 12D  &13A & 14A  \\

\midrule
-3	&2	&2	&2	&2	&2	&2	&2	&2	&2	&2	&2	&2	&2	&2	&2	&2	&2	&2	&2	&2	&2 &2  &2  &2	\\
0       & 248& -8& 14& 5& -4& 8& 0& -2& 4& -2& 1& 3& 0& 0& 5& -4& 2& 2& 2& 2& -1& 0& 1& -1 \\
4      & 54000& 240& -54& 54& 0& 16& 0& 0& 0& -6& 6& 2& 0& 0& 0& 0& 0& 0& -2&  -2& -2& 0& -2& 2 \\
5  &  -171990& 42& 0& 54& -54& 42& -22& 10& -6& 0& 6& 0& -6& 2& 0& 0& 0& 2&  0& 0& 6& 2& 0& 0 \\
8  & 3414528& -1536& 0& -108& 108& 0& 0& 28& -12& 0& 12& 12& 0& 0& 0& 0&  0& 4& 0& 0& 0& 0& 0& 4 \\
9 & -8192000& 0& -128& 16& 160& 0& 0& 0& 0& 0& 0& 16& 0& 0& 16& -2& -8&  0& 0& 0& 0& 0& 2& 0 \\
12 & 88660992& 7168& 336& 12& -312& 0& 0& -8& -8& 16& 4& 0& 0& 0& 12& 12& -6& 8& 0& 0& 0& 0& 4& 0 \\
13 & -183902292& -84& 378& -378& 0& -84& 108& -42& 0& -6& 6& 0& -4& -12&  0& 0& 0& 6& -6& -6& 6& 0& -2& 0 \\
16 & 1417878000& -26640& -702& 702& 0& 16& 0& 0& 0& 18& -18& -14& 0& 0& 0&  0& 0& 0& -2& -2& -2& 0& 6& 2 \\
17 & -2687827968& 0& 0& 864& -864& 0& 0& 32& 0& 0& 0& 0& 0& 0& 0& 0& 0& 0&  0& 0& 0& 0& 6& 0 \\
20 & 16555069440& 86016& 0& -1512& 1512& 0& 0& -60& 24& 0& -24& 0& 0& 0&  0& 0& 0& -4& 0& 0& 0& 0& 0& 0 \\
21 & 29466050250& 310& -1794& -12& 1770& 310& -330& 0& 10& -2& 4& -14& 22& -10& -12& -12& 6& 0& -2& -2& 4& -6& 0& 2 \\
24 & 154184576000& -250880& 2912& -4& -2920& 0& 0& 0& 40& -32& 4& 0& 0& 0& -4& -4& -4& 0& 0& 0& 0& 0& 0& 0 \\
25 & -261761531904& 0& 3456& -3456& 0& 0& 0& 96& 0& 0& 0& 16& 0& 0& 0& 0&  0& 0& 0& 0& 0& 0& -2& 0 \\
28 & 1208278536192& 675840& -5616& 5616& 0& 0& 0& 192& 0& -48& 48& 28& 0&  0& 0& 0& 0& 0& 0& 0& 0& 0& 0& 4 \\
29 & -1976452670250& -810& 0& 6480& -6480& -810& 790& 0& 0& 0& 0& -30& 6&  30& 0& 0& 0& 0& 0& 0& 0& -8& -6& 2 \\
32 & 8251985424384& -1708032& 0& -10044& 10044& 0& 0& -116& -60& 0& 60&  -12& 0& 0& 0& 0& 0& -12& 0& 0& 0& 0& 0& -4 \\
\bottomrule
\end{tabular}
\end{small}
\end{sidewaystable}

\begin{sidewaystable}
\begin{small}
\centering
\caption{Coefficient of $q^n$ in the McKay--Thompson series ${\cal F}_{3,[g]}(\tau)$ of Thompson moonshine, part two.}\label{miseries}\smallskip
\begin{tabular}{c@{ }@{\;}r@{ }r@{ }r@{ }r@{ }r@{ }r@{ }r@{ }r@{ }r@{ }r@{ }r@{ }r@{ }r@{ }r@{ }r@{ }r@{ }r@{ }r@{ }r@{ }r@{ }r@{ }r@{ }r@{ }r@{ }r@{ }r}\toprule
$n~ \backslash~ [g]$	&15A&15B&18A&18B&19A&20A&21A	&24A	&24B	&24C	&24D	&27A	&27B&27C&28A&30A	&30B	&31A	&31B	&36A  & 36B	&36C  & 39A &39B &   \\
\midrule
-3	&2	&2	&2	&2	&2	&2	&2	&2	&2	&2	&2	&2	&2	&2	&2	&2	&2	&2	&2	&2	&2	&2 &2 &2 \\
0     & 1& 1& 1& -2& 1& 0& 0& 0& 0& 0& 0& 2& -1& -1& 1& -1& -1& 0& 0& -1& -1&  -1& 1& 1 \\
4     &  0& 0& 0& 0& 2& 0& 2& 0& 0& 0& 0& 0& 0& 0& 2& 0& 0& -2& -2& 4& -2& -2&  -2& -2 \\
5   &  1& 1& 0& 0& -2& -2& 0& 0& 0& 2& 2& 0& 0& 0& 0& -1& -1& -2& -2& 0& 0&  0& 0& 0 \\
8   &  -2& -2& 0& 0& 0& 0& 0& 0& 0& 0& 0& 0& 0& 0& 0& -2& -2& 2& 2& 0& 0& 0&  0& 0 \\
9   &  0& 0& 0& 0& 2& 0& -2& 0& 0& 0& 0& -2& 1& 1& 0& 0& 0& -2& -2& 0& 0& 0&  2& 2 \\
12  & -2& -2& 4& -2& 0& 0& 0& 0& 0& 0& 0& 0& 0& 0& 0& 2& 2& 0& 0& 0& 0& 0&  -2& -2 \\
13 & 0& 0& 0& 0& 0& -2& 0& 2& 2& 0& 0& 0& 0& 0& 0& 0& 0& 0& 0& 0& 0& 0& 1& 1 \\
16 & 0& 0& 0& 0& -2& 0& -2& 0& 0& 0& 0& 0& 0& 0& 2& 0& 0& 0& 0& 4& -2& -2&  0& 0 \\
17 & -4& -4& 0& 0& 2& 0& 0& 0& 0& 0& 0& 0& 0& 0& 0& 0& 0& 0& 0& 0& 0& 0& 0& 0 \\
20 & -3& -3& 0& 0& 4& 0& 0& 0& 0& 0& 0& 0& 0& 0& 0& -1& -1& 2& 2& 0& 0& 0&  0& 0 \\
21 & 0& 0& 4& -2& 0& 0& -2& -2& -2& 2& 2& 0& 0& 0& 2& 0& 0& 0& 0& 4& 4& 4&  0& 0 \\
24 & 0& 0& 4& 4& 4& 0& 0& 0& 0& 0& 0& -4& -4& -4& 0& 0& 0& 0& 0& 0& 0& 0&  0& 0 \\
25 & 0& 0& 0& 0& 4& 0& -2& 0& 0& 0& 0& 0& 0& 0& 0& 0& 0& 0& 0& 0& 0& 0& -2& -2 \\
28 & 0& 0& 0& 0& 0& 0& -2& 0& 0& 0& 0& 0& 0& 0& 0& 0& 0& -2& -2& 0& 0& 0&  0& 0 \\
29 & 0& 0& 0& 0& 0& 0& 0& 0& 0& 0& 0& 0& 0& 0& 2& 0& 0& 0& 0& 0& 0& 0& 0& 0 \\
32 & 4& 4& 0& 0& 0& 0& 0& 0& 0& 0& 0& 0& 0& 0& 0& 0& 0& -2& -2& 0& 0& 0& 0& 0 \\
\bottomrule
\end{tabular}
\end{small}
\end{sidewaystable}
\begin{sidewaystable}
\begin{center}
\caption{\label{thchartabone} {Character table of the Thompson group, part one. The notation follows that of the ATLAS \cite{atlas} and the GAP software package \cite{GAP4} with $bN=(-1+ i\sqrt{N})/2$ for $N \equiv -1 ~{\rm mod}~4$ and $iN=i \sqrt{N}$, where $i= \sqrt{-1}$ and $A=1+4b3$, $B=2+8 b3$, $C=b15$, $D=-i3$,$E=i6$, $F= 1+3 b3$, $G=b31$, $H=2b3$ and $I=b39$. An overline indicates complex conjugation.}}
\smallskip
\begin{small}
\begin{tabular}{c@{ }r@{ }r@{ }r@{ }r@{ }r@{ }r@{ }r@{ }r@{ }r@{ }r@{ }r@{ }r@{ }r@{ }r@{ }r@{ }r@{ }r@{ }r@{ }r@{ }r@{ }r@{ }r@{ }r@{ }r} \toprule
$[g]$ & 1A  & 2A  & 3A  & 3B  & 3C  & 4A  & 4B  & 5A  & 6A  & 6B  & 6C  & 7A  & 8A  & 8B  & 9A  & 9B  & 9C  & 10A  & 12A  & 12B  & 12C  & 12D  & 13A  & 14A  \\ 
	 \midrule$\chi_{1}$ & $1$  & $1$  & $1$  & $1$  & $1$  & $1$  & $1$  & $1$  & $1$  & $1$  & $1$  & $1$  & $1$  & $1$  & $1$  & $1$  & $1$  & $1$  & $1$  & $1$  & $1$  & $1$  & $1$  & $1$ \\
		 $\chi_{2}$ & $248$  & -$8$  & $14$  & $5$  & -$4$  & $8$  & $0$  & -$2$  & $4$  & -$2$  & $1$  & $3$  & $0$  & $0$  & $5$  & -$4$  & $2$  & $2$  & $2$  & $2$  & -$1$  & $0$  & $1$  & -$1$ \\
		 $\chi_{3}$ & $4123$  & $27$  & $64$  & -$8$  & $1$  & $27$  & -$5$  & -$2$  & $9$  & $0$  & $0$  & $7$  & $3$  & -$1$  & -$8$  & $1$  & $4$  & $2$  & $0$  & $0$  & $0$  & $1$  & $2$  & -$1$ \\
		 $\chi_{4}$ & $27000$  & $120$  & -$27$  & $27$  & $0$  & $8$  & $0$  & $0$  & $0$  & -$3$  & $3$  & $1$  & $0$  & $0$  & $0$  & $0$  & $0$  & $0$  & $A$  & $\overline{A}$ & -$1$  & $0$  & -$1$  & $1$ \\
		 $\chi_{5}$ & $27000$  & $120$  & -$27$  & $27$  & $0$  & $8$  & $0$  & $0$  & $0$  & -$3$  & $3$  & $1$  & $0$  & $0$  & $0$  & $0$  & $0$  & $0$  & $\overline{A}$  & $A$  & -$1$  & $0$  & -$1$  & $1$ \\
		 $\chi_{6}$ & $30628$  & -$92$  & $91$  & $10$  & $10$  & $36$  & $4$  & $3$  & $10$  & -$5$  & -$2$  & $3$  & -$4$  & $0$  & $10$  & $10$  & $1$  & $3$  & $3$  & $3$  & $0$  & -$2$  & $0$  & -$1$ \\
		 $\chi_{7}$ & $30875$  & $155$  & $104$  & $14$  & $5$  & $27$  & -$5$  & $0$  & $5$  & $8$  & $2$  & $5$  & $3$  & -$1$  & $14$  & $5$  & $2$  & $0$  & $0$  & $0$  & $0$  & $1$  & $0$  & $1$ \\
		 $\chi_{8}$ & $61256$  & $72$  & $182$  & $20$  & $20$  & $56$  & $0$  & $6$  & $12$  & $6$  & $0$  & $6$  & $0$  & $0$  & -$7$  & -$7$  & $2$  & $2$  & $2$  & $2$  & $2$  & $0$  & $0$  & $2$ \\
		 $\chi_{9}$ & $85995$  & -$21$  & $0$  & -$27$  & $27$  & -$21$  & $11$  & -$5$  & $3$  & $0$  & -$3$  & $0$  & $3$  & -$1$  & $0$  & $0$  & $0$  & -$1$  & $0$  & $0$  & -$3$  & -$1$  & $0$  & $0$ \\
		 $\chi_{10}$ & $85995$  & -$21$  & $0$  & -$27$  & $27$  & -$21$  & $11$  & -$5$  & $3$  & $0$  & -$3$  & $0$  & $3$  & -$1$  & $0$  & $0$  & $0$  & -$1$  & $0$  & $0$  & -$3$  & -$1$  & $0$  & $0$ \\
		 $\chi_{11}$ & $147250$  & $50$  & $181$  & -$8$  & -$35$  & $34$  & $10$  & $0$  & $5$  & $5$  & -$4$  & $5$  & $2$  & -$2$  & $19$  & -$8$  & $1$  & $0$  & $1$  & $1$  & -$2$  & $1$  & -$1$  & $1$ \\
		 $\chi_{12}$ & $767637$  & $405$  & $0$  & $0$  & $0$  & -$27$  & -$3$  & $12$  & $0$  & $0$  & $0$  & $3$  & -$3$  & -$3$  & $0$  & $0$  & $0$  & $0$  & $0$  & $0$  & $0$  & $0$  & $0$  & -$1$ \\
		 $\chi_{13}$ & $767637$  & $405$  & $0$  & $0$  & $0$  & -$27$  & -$3$  & $12$  & $0$  & $0$  & $0$  & $3$  & -$3$  & -$3$  & $0$  & $0$  & $0$  & $0$  & $0$  & $0$  & $0$  & $0$  & $0$  & -$1$ \\
		 $\chi_{14}$ & $779247$  & -$273$  & -$189$  & -$54$  & $0$  & $63$  & -$9$  & -$3$  & $0$  & $3$  & $6$  & $0$  & -$1$  & $3$  & $0$  & $0$  & $0$  & -$3$  & $3$  & $3$  & $0$  & $0$  & $1$  & $0$ \\
		 $\chi_{15}$ & $779247$  & -$273$  & -$189$  & -$54$  & $0$  & $63$  & -$9$  & -$3$  & $0$  & $3$  & $6$  & $0$  & -$1$  & $3$  & $0$  & $0$  & $0$  & -$3$  & $3$  & $3$  & $0$  & $0$  & $1$  & $0$ \\
		 $\chi_{16}$ & $957125$  & -$315$  & $650$  & -$52$  & -$25$  & $133$  & $5$  & $0$  & $15$  & -$6$  & $0$  & $8$  & -$3$  & $1$  & -$25$  & $2$  & $2$  & $0$  & -$2$  & -$2$  & -$2$  & -$1$  & $0$  & $0$ \\
		 $\chi_{17}$ & $1707264$  & -$768$  & $0$  & -$54$  & $54$  & $0$  & $0$  & $14$  & -$6$  & $0$  & $6$  & $6$  & $0$  & $0$  & $0$  & $0$  & $0$  & $2$  & $0$  & $0$  & $0$  & $0$  & $0$  & $2$ \\
		 $\chi_{18}$ & $1707264$  & -$768$  & $0$  & -$54$  & $54$  & $0$  & $0$  & $14$  & -$6$  & $0$  & $6$  & $6$  & $0$  & $0$  & $0$  & $0$  & $0$  & $2$  & $0$  & $0$  & $0$  & $0$  & $0$  & $2$ \\
		 $\chi_{19}$ & $2450240$  & $832$  & $260$  & $71$  & $44$  & $64$  & $0$  & -$10$  & $4$  & $4$  & -$5$  & -$5$  & $0$  & $0$  & $17$  & -$10$  & -$1$  & $2$  & $4$  & $4$  & $1$  & $0$  & $0$  & -$1$ \\
		 $\chi_{20}$ & $2572752$  & -$1072$  & $624$  & $111$  & $84$  & $48$  & $0$  & $2$  & -$4$  & -$16$  & -$1$  & $7$  & $0$  & $0$  & $30$  & $3$  & $3$  & -$2$  & $0$  & $0$  & $3$  & $0$  & $0$  & -$1$ \\
		 $\chi_{21}$ & $3376737$  & $609$  & $819$  & $9$  & $9$  & $161$  & $1$  & -$13$  & $9$  & $3$  & -$3$  & $0$  & $1$  & $1$  & $9$  & $9$  & $0$  & -$1$  & -$1$  & -$1$  & -$1$  & $1$  & $0$  & $0$ \\
		 $\chi_{22}$ & $4096000$  & $0$  & $64$  & -$8$  & -$80$  & $0$  & $0$  & $0$  & $0$  & $0$  & $0$  & -$8$  & $0$  & $0$  & -$8$  & $1$  & $4$  & $0$  & $0$  & $0$  & $0$  & $0$  & -$1$  & $0$ \\
		 $\chi_{23}$ & $4096000$  & $0$  & $64$  & -$8$  & -$80$  & $0$  & $0$  & $0$  & $0$  & $0$  & $0$  & -$8$  & $0$  & $0$  & -$8$  & $1$  & $4$  & $0$  & $0$  & $0$  & $0$  & $0$  & -$1$  & $0$ \\
		 $\chi_{24}$ & $4123000$  & $120$  & $118$  & $19$  & -$80$  & $8$  & $0$  & $0$  & $0$  & $6$  & $3$  & -$7$  & $0$  & $0$  & $19$  & $1$  & $4$  & $0$  & $2$  & $2$  & -$1$  & $0$  & -$2$  & $1$ \\
		 \bottomrule
	\end{tabular}
	\end{small}
	\end{center}
		\end{sidewaystable}
	
\begin{sidewaystable}
\begin{center}
\caption{Character table of the Thompson group, part two.}\label{thchartabtwo}
\smallskip
\begin{small}
\begin{tabular}{c@{ }r@{ }r@{ }r@{ }r@{ }r@{ }r@{ }r@{ }r@{ }r@{ }r@{ }r@{ }r@{ }r@{ }r@{ }r@{ }r@{ }r@{ }r@{ }r@{ }r@{ }r@{ }r@{ }r@{ }r} \toprule
$[g]$ & 15A  & 15B  & 18A  & 18B  & 19A  & 20A  & 21A  & 24A  & 24B  & 24C  & 24D  & 27A  & 27B  & 27C  & 28A  & 30A  & 30B  & 31A  & 31B  & 36A  & 36B  & 36C  & 39A  & 39B  \\ 
	 \midrule$\chi_{1}$ & $1$  & $1$  & $1$  & $1$  & $1$  & $1$  & $1$  & $1$  & $1$  & $1$  & $1$  & $1$  & $1$  & $1$  & $1$  & $1$  & $1$  & $1$  & $1$  & $1$  & $1$  & $1$  & $1$  & $1$ \\
		 $\chi_{2}$ & $1$  & $1$  & $1$  & -$2$  & $1$  & $0$  & $0$  & $0$  & $0$  & $0$  & $0$  & $2$  & -$1$  & -$1$  & $1$  & -$1$  & -$1$  & $0$  & $0$  & -$1$  & -$1$  & -$1$  & $1$  & $1$ \\
		 $\chi_{3}$ & $1$  & $1$  & $0$  & $0$  & $0$  & $0$  & $1$  & $0$  & $0$  & -$1$  & -$1$  & -$2$  & $1$  & $1$  & -$1$  & -$1$  & -$1$  & $0$  & $0$  & $0$  & $0$  & $0$  & -$1$  & -$1$ \\
		 $\chi_{4}$ & $0$  & $0$  & $0$  & $0$  & $1$  & $0$  & $1$  & -$D$  & $\overline{D}$  & $0$  & $0$  & $0$  & $0$  & $0$  & $1$  & $0$  & $0$  & -$1$  & -$1$  & $2$  & $H$  & $\overline{H}$  & -$1$  & -$1$ \\
		 $\chi_{5}$ & $0$  & $0$  & $0$  & $0$  & $1$  & $0$  & $1$  & $\overline{D}$  & $D$  & $0$  & $0$  & $0$  & $0$  & $0$  & $1$  & $0$  & $0$  & -$1$  & -$1$  & $2$  & $\overline{H}$  & $H$  & -$1$  & -$1$ \\
		 $\chi_{6}$ & $0$  & $0$  & -$2$  & $1$  & $0$  & -$1$  & $0$  & -$1$  & -$1$  & $0$  & $0$  & $1$  & $1$  & $1$  & $1$  & $0$  & $0$  & $0$  & $0$  & $0$  & $0$  & $0$  & $0$  & $0$ \\
		 $\chi_{7}$ & $0$  & $0$  & $2$  & $2$  & $0$  & $0$  & -$1$  & $0$  & $0$  & -$1$  & -$1$  & $2$  & -$1$  & -$1$  & -$1$  & $0$  & $0$  & -$1$  & -$1$  & $0$  & $0$  & $0$  & $0$  & $0$ \\
		 $\chi_{8}$ & $0$  & $0$  & -$3$  & $0$  & $0$  & $0$  & $0$  & $0$  & $0$  & $0$  & $0$  & -$1$  & -$1$  & -$1$  & $0$  & $2$  & $2$  & $0$  & $0$  & -$1$  & -$1$  & -$1$  & $0$  & $0$ \\
		 $\chi_{9}$ & $C$  & $\overline{C}$ & $0$  & $0$  & $1$  & $1$  & $0$  & $0$  & $0$  & -$1$  & -$1$  & $0$  & $0$  & $0$  & $0$  & -$C$  & -$ \overline{C}$  & $1$  & $1$  & $0$  & $0$  & $0$  & $0$  & $0$ \\
		 $\chi_{10}$ & $\overline{C}$  & $C$  & $0$  & $0$  & $1$  & $1$  & $0$  & $0$  & $0$  & -$1$  & -$1$  & $0$  & $0$  & $0$  & $0$  & -$\overline{C}$  & -$C$  & $1$  & $1$  & $0$  & $0$  & $0$  & $0$  & $0$ \\
		 $\chi_{11}$ & $0$  & $0$  & -$1$  & -$1$  & $0$  & $0$  & -$1$  & -$1$  & -$1$  & $1$  & $1$  & $1$  & $1$  & $1$  & -$1$  & $0$  & $0$  & $0$  & $0$  & $1$  & $1$  & $1$  & -$1$  & -$1$ \\
		 $\chi_{12}$ & $0$  & $0$  & $0$  & $0$  & -$1$  & $2$  & $0$  & $0$  & $0$  & $0$  & $0$  & $0$  & $0$  & $0$  & $1$  & $0$  & $0$  & $G$  & $\overline{G}$  & $0$  & $0$  & $0$  & $0$  & $0$ \\
		 $\chi_{13}$ & $0$  & $0$  & $0$  & $0$  & -$1$  & $2$  & $0$  & $0$  & $0$  & $0$  & $0$  & $0$  & $0$  & $0$  & $1$  & $0$  & $0$  & $\overline{G}$  & $G$  & $0$  & $0$  & $0$  & $0$  & $0$ \\
		 $\chi_{14}$ & $0$  & $0$  & $0$  & $0$  & $0$  & $1$  & $0$  & -$1$  & -$1$  & $0$  & $0$  & $0$  & $0$  & $0$  & $0$  & $0$  & $0$  & $0$  & $0$  & $0$  & $0$  & $0$  & $I$  & $\overline{I}$ \\
		 $\chi_{15}$ & $0$  & $0$  & $0$  & $0$  & $0$  & $1$  & $0$  & -$1$  & -$1$  & $0$  & $0$  & $0$  & $0$  & $0$  & $0$  & $0$  & $0$  & $0$  & $0$  & $0$  & $0$  & $0$  & $\overline{I}$  & $I$ \\
		 $\chi_{16}$ & $0$  & $0$  & $3$  & $0$  & $0$  & $0$  & -$1$  & $0$  & $0$  & $1$  & $1$  & -$1$  & -$1$  & -$1$  & $0$  & $0$  & $0$  & $0$  & $0$  & $1$  & $1$  & $1$  & $0$  & $0$ \\
		 $\chi_{17}$ & -$1$  & -$1$  & $0$  & $0$  & $0$  & $0$  & $0$  & $0$  & $0$  & $E$  & $\overline{E}$  & $0$  & $0$  & $0$  & $0$  & -$1$  & -$1$  & $1$  & $1$  & $0$  & $0$  & $0$  & $0$  & $0$ \\
		 $\chi_{18}$ & -$1$  & -$1$  & $0$  & $0$  & $0$  & $0$  & $0$  & $0$  & $0$  & $\overline{E}$  & $E$  & $0$  & $0$  & $0$  & $0$  & -$1$  & -$1$  & $1$  & $1$  & $0$  & $0$  & $0$  & $0$  & $0$ \\
		 $\chi_{19}$ & -$1$  & -$1$  & $1$  & $1$  & $0$  & $0$  & $1$  & $0$  & $0$  & $0$  & $0$  & -$1$  & -$1$  & -$1$  & $1$  & -$1$  & -$1$  & $0$  & $0$  & $1$  & $1$  & $1$  & $0$  & $0$ \\
		 $\chi_{20}$ & -$1$  & -$1$  & $2$  & -$1$  & $0$  & $0$  & $1$  & $0$  & $0$  & $0$  & $0$  & $0$  & $0$  & $0$  & -$1$  & $1$  & $1$  & $0$  & $0$  & $0$  & $0$  & $0$  & $0$  & $0$ \\
		 $\chi_{21}$ & -$1$  & -$1$  & -$3$  & $0$  & $0$  & $1$  & $0$  & $1$  & $1$  & $1$  & $1$  & $0$  & $0$  & $0$  & $0$  & -$1$  & -$1$  & $0$  & $0$  & -$1$  & -$1$  & -$1$  & $0$  & $0$ \\
		 $\chi_{22}$ & $0$  & $0$  & $0$  & $0$  & -$1$  & $0$  & $1$  & $0$  & $0$  & $0$  & $0$  & $1$  & $F$  & $\overline{F}$  & $0$  & $0$  & $0$  & $1$  & $1$  & $0$  & $0$  & $0$  & -$1$  & -$1$ \\
		 $\chi_{23}$ & $0$  & $0$  & $0$  & $0$  & -$1$  & $0$  & $1$  & $0$  & $0$  & $0$  & $0$  & $1$  & $\overline{F}$  & $F$  & $0$  & $0$  & $0$  & $1$  & $1$  & $0$  & $0$  & $0$  & -$1$  & -$1$ \\
		 $\chi_{24}$ & $0$  & $0$  & $3$  & $0$  & $0$  & $0$  & -$1$  & $0$  & $0$  & $0$  & $0$  & -$2$  & $1$  & $1$  & $1$  & $0$  & $0$  & $0$  & $0$  & -$1$  & -$1$  & -$1$  & $1$  & $1$ \\
		 \bottomrule
	\end{tabular}
	\end{small}
	\end{center}
	\end{sidewaystable}
	
\begin{sidewaystable}
\begin{center}
\caption{Character table of the Thompson group, part three. }\label{thchartabthree}
\smallskip
\begin{small}
\begin{tabular}{c@{ }r@{ }r@{ }r@{ }r@{ }r@{ }r@{ }r@{ }r@{ }r@{ }r@{ }r@{ }r@{ }r@{ }r@{ }r@{ }r@{ }r@{ }r@{ }r@{ }r@{ }r@{ }r@{ }r@{ }r} \toprule
$[g]$ & 1A  & 2A  & 3A  & 3B  & 3C  & 4A  & 4B  & 5A  & 6A  & 6B  & 6C  & 7A  & 8A  & 8B  & 9A  & 9B  & 9C  & 10A  & 12A  & 12B  & 12C  & 12D  & 13A  & 14A  \\ 
	 \midrule$\chi_{25}$ & $4881384$  & $1512$  & $729$  & $0$  & $0$  & $72$  & $24$  & $9$  & $0$  & $9$  & $0$  & $4$  & $8$  & $0$  & $0$  & $0$  & $0$  & -$3$  & -$3$  & -$3$  & $0$  & $0$  & $1$  & $0$ \\
		 $\chi_{26}$ & $4936750$  & -$210$  & $637$  & -$38$  & -$65$  & $126$  & -$10$  & $0$  & $15$  & -$3$  & $6$  & $0$  & -$2$  & $2$  & $16$  & -$11$  & -$2$  & $0$  & -$3$  & -$3$  & $0$  & -$1$  & $0$  & $0$ \\
		 $\chi_{27}$ & $6669000$  & -$1080$  & -$351$  & $108$  & $0$  & $56$  & $0$  & $0$  & $0$  & $9$  & $0$  & $2$  & $0$  & $0$  & $0$  & $0$  & $0$  & $0$  & $A$  & $\overline{A}$  & $2$  & $0$  & $0$  & -$2$ \\
		 $\chi_{28}$ & $6669000$  & -$1080$  & -$351$  & $108$  & $0$  & $56$  & $0$  & $0$  & $0$  & $9$  & $0$  & $2$  & $0$  & $0$  & $0$  & $0$  & $0$  & $0$  & $\overline{A}$  & $A$  & $2$  & $0$  & $0$  & -$2$ \\
		 $\chi_{29}$ & $6696000$  & -$960$  & -$378$  & $135$  & $0$  & $64$  & $0$  & $0$  & $0$  & $6$  & $3$  & $3$  & $0$  & $0$  & $0$  & $0$  & $0$  & $0$  & $B$  & $\overline{B}$  & $1$  & $0$  & -$1$  & -$1$ \\
		 $\chi_{30}$ & $6696000$  & -$960$  & -$378$  & $135$  & $0$  & $64$  & $0$  & $0$  & $0$  & $6$  & $3$  & $3$  & $0$  & $0$  & $0$  & $0$  & $0$  & $0$  & $\overline{B}$  & $B$  & $1$  & $0$  & -$1$  & -$1$ \\
		 $\chi_{31}$ & $10822875$  & -$805$  & $924$  & $141$  & -$75$  & $91$  & -$5$  & $0$  & $5$  & -$4$  & $5$  & $0$  & $3$  & -$1$  & -$21$  & $6$  & -$3$  & $0$  & $4$  & $4$  & $1$  & $1$  & -$2$  & $0$ \\
		 $\chi_{32}$ & $11577384$  & $552$  & $351$  & $135$  & $0$  & -$120$  & $24$  & $9$  & $0$  & $15$  & $3$  & $7$  & -$8$  & $0$  & $0$  & $0$  & $0$  & -$3$  & $3$  & $3$  & -$3$  & $0$  & $0$  & -$1$ \\
		 $\chi_{33}$ & $16539120$  & $2544$  & $0$  & $297$  & -$54$  & $48$  & $16$  & -$5$  & -$6$  & $0$  & -$3$  & $3$  & $0$  & $0$  & $0$  & $0$  & $0$  & -$1$  & $0$  & $0$  & $3$  & -$2$  & $0$  & $3$ \\
		 $\chi_{34}$ & $18154500$  & $1540$  & -$273$  & $213$  & -$30$  & -$28$  & $20$  & $0$  & $10$  & -$17$  & $1$  & $0$  & -$4$  & $0$  & -$3$  & -$3$  & -$3$  & $0$  & -$1$  & -$1$  & -$1$  & $2$  & $0$  & $0$ \\
		 $\chi_{35}$ & $21326760$  & $168$  & $0$  & -$135$  & -$108$  & -$168$  & $0$  & $10$  & $12$  & $0$  & -$3$  & $0$  & $0$  & $0$  & $0$  & $0$  & $0$  & -$2$  & $0$  & $0$  & $3$  & $0$  & $0$  & $0$ \\
		 $\chi_{36}$ & $21326760$  & $168$  & $0$  & -$135$  & -$108$  & -$168$  & $0$  & $10$  & $12$  & $0$  & -$3$  & $0$  & $0$  & $0$  & $0$  & $0$  & $0$  & -$2$  & $0$  & $0$  & $3$  & $0$  & $0$  & $0$ \\
		 $\chi_{37}$ & $28861000$  & $840$  & $1078$  & -$110$  & $160$  & $56$  & $0$  & $0$  & $0$  & $6$  & -$6$  & $0$  & $0$  & $0$  & -$29$  & -$2$  & -$2$  & $0$  & $2$  & $2$  & $2$  & $0$  & -$1$  & $0$ \\
		 $\chi_{38}$ & $30507008$  & $0$  & $896$  & -$184$  & $32$  & $0$  & $0$  & $8$  & $0$  & $0$  & $0$  & $0$  & $0$  & $0$  & $32$  & $5$  & -$4$  & $0$  & $0$  & $0$  & $0$  & $0$  & -$1$  & $0$ \\
		 $\chi_{39}$ & $40199250$  & $3410$  & -$78$  & $3$  & $165$  & -$62$  & $10$  & $0$  & $5$  & $2$  & -$1$  & -$7$  & -$6$  & $2$  & $3$  & $3$  & $3$  & $0$  & -$2$  & -$2$  & $1$  & $1$  & $0$  & $1$ \\
		 $\chi_{40}$ & $44330496$  & $3584$  & $168$  & $6$  & -$156$  & $0$  & $0$  & -$4$  & -$4$  & $8$  & $2$  & $0$  & $0$  & $0$  & $6$  & $6$  & -$3$  & $4$  & $0$  & $0$  & $0$  & $0$  & $2$  & $0$ \\
		 $\chi_{41}$ & $51684750$  & $2190$  & $0$  & $108$  & $135$  & -$162$  & -$10$  & $0$  & $15$  & $0$  & $12$  & -$9$  & $6$  & -$2$  & $0$  & $0$  & $0$  & $0$  & $0$  & $0$  & $0$  & -$1$  & $0$  & -$1$ \\
		 $\chi_{42}$ & $72925515$  & -$2997$  & $0$  & $0$  & $0$  & $27$  & $51$  & $15$  & $0$  & $0$  & $0$  & -$9$  & $3$  & $3$  & $0$  & $0$  & $0$  & $3$  & $0$  & $0$  & $0$  & $0$  & $0$  & -$1$ \\
		 $\chi_{43}$ & $76271625$  & -$2295$  & $729$  & $0$  & $0$  & $153$  & -$15$  & $0$  & $0$  & $9$  & $0$  & -$11$  & -$7$  & -$3$  & $0$  & $0$  & $0$  & $0$  & -$3$  & -$3$  & $0$  & $0$  & $1$  & $1$ \\
		 $\chi_{44}$ & $77376000$  & $2560$  & $1560$  & -$60$  & -$60$  & $0$  & $0$  & $0$  & -$20$  & -$8$  & $4$  & $2$  & $0$  & $0$  & -$6$  & -$6$  & $3$  & $0$  & $0$  & $0$  & $0$  & $0$  & $0$  & -$2$ \\
		 $\chi_{45}$ & $81153009$  & -$783$  & -$729$  & $0$  & $0$  & $225$  & $9$  & $9$  & $0$  & -$9$  & $0$  & -$7$  & $1$  & -$3$  & $0$  & $0$  & $0$  & -$3$  & $3$  & $3$  & $0$  & $0$  & $2$  & $1$ \\
		 $\chi_{46}$ & $91171899$  & $315$  & $0$  & $243$  & $0$  & -$21$  & -$45$  & $24$  & $0$  & $0$  & -$9$  & $0$  & $3$  & $3$  & $0$  & $0$  & $0$  & $0$  & $0$  & $0$  & -$3$  & $0$  & $0$  & $0$ \\
		 $\chi_{47}$ & $111321000$  & $3240$  & -$1728$  & -$216$  & $0$  & $216$  & $0$  & $0$  & $0$  & $0$  & $0$  & $7$  & $0$  & $0$  & $0$  & $0$  & $0$  & $0$  & $0$  & $0$  & $0$  & $0$  & -$2$  & -$1$ \\
		 $\chi_{48}$ & $190373976$  & -$3240$  & $0$  & $0$  & $0$  & -$216$  & $0$  & -$24$  & $0$  & $0$  & $0$  & $9$  & $0$  & $0$  & $0$  & $0$  & $0$  & $0$  & $0$  & $0$  & $0$  & $0$  & $0$  & $1$ \\
		 \bottomrule
	\end{tabular}
	\end{small}
	\end{center}
	\end{sidewaystable}
	
\begin{sidewaystable}
\begin{center}
\caption{Character table of the Thompson group, part four. }\label{thchartabfour}
\smallskip
\begin{small}
\begin{tabular}{c@{ }r@{ }r@{ }r@{ }r@{ }r@{ }r@{ }r@{ }r@{ }r@{ }r@{ }r@{ }r@{ }r@{ }r@{ }r@{ }r@{ }r@{ }r@{ }r@{ }r@{ }r@{ }r@{ }r@{ }r} \toprule
$[g]$ & 15A  & 15B  & 18A  & 18B  & 19A  & 20A  & 21A  & 24A  & 24B  & 24C  & 24D  & 27A  & 27B  & 27C  & 28A  & 30A  & 30B  & 31A  & 31B  & 36A  & 36B  & 36C  & 39A  & 39B  \\ 
	 \midrule$\chi_{25}$ & $0$  & $0$  & $0$  & $0$  & -$1$  & -$1$  & $1$  & -$1$  & -$1$  & $0$  & $0$  & $0$  & $0$  & $0$  & $2$  & $0$  & $0$  & $0$  & $0$  & $0$  & $0$  & $0$  & $1$  & $1$ \\
		 $\chi_{26}$ & $0$  & $0$  & $0$  & $0$  & -$1$  & $0$  & $0$  & $1$  & $1$  & -$1$  & -$1$  & $1$  & $1$  & $1$  & $0$  & $0$  & $0$  & $0$  & $0$  & $0$  & $0$  & $0$  & $0$  & $0$ \\
		 $\chi_{27}$ & $0$  & $0$  & $0$  & $0$  & $0$  & $0$  & -$1$  & $\overline{D}$  & $D$  & $0$  & $0$  & $0$  & $0$  & $0$  & $0$  & $0$  & $0$  & $1$  & $1$  & $2$  & $H$  & $\overline{H}$  & $0$  & $0$ \\
		 $\chi_{28}$ & $0$  & $0$  & $0$  & $0$  & $0$  & $0$  & -$1$  & $D$  & $\overline{D}$  & $0$  & $0$  & $0$  & $0$  & $0$  & $0$  & $0$  & $0$  & $1$  & $1$  & $2$  & $\overline{H}$  & $H$  & $0$  & $0$ \\
		 $\chi_{29}$ & $0$  & $0$  & $0$  & $0$  & $1$  & $0$  & $0$  & $0$  & $0$  & $0$  & $0$  & $0$  & $0$  & $0$  & $1$  & $0$  & $0$  & $0$  & $0$  & -$2$  & -$H$  & -$\overline{H}$  & -$1$  & -$1$ \\
		 $\chi_{30}$ & $0$  & $0$  & $0$  & $0$  & $1$  & $0$  & $0$  & $0$  & $0$  & $0$  & $0$  & $0$  & $0$  & $0$  & $1$  & $0$  & $0$  & $0$  & $0$  & -$2$  & -$\overline{H}$  & -$H$  & -$1$  & -$1$ \\
		 $\chi_{31}$ & $0$  & $0$  & -$1$  & -$1$  & $0$  & $0$  & $0$  & $0$  & $0$  & -$1$  & -$1$  & $0$  & $0$  & $0$  & $0$  & $0$  & $0$  & $0$  & $0$  & $1$  & $1$  & $1$  & $1$  & $1$ \\
		 $\chi_{32}$ & $0$  & $0$  & $0$  & $0$  & $0$  & -$1$  & $1$  & $1$  & $1$  & $0$  & $0$  & $0$  & $0$  & $0$  & -$1$  & $0$  & $0$  & $0$  & $0$  & $0$  & $0$  & $0$  & $0$  & $0$ \\
		 $\chi_{33}$ & $1$  & $1$  & $0$  & $0$  & $0$  & $1$  & $0$  & $0$  & $0$  & $0$  & $0$  & $0$  & $0$  & $0$  & -$1$  & -$1$  & -$1$  & $0$  & $0$  & $0$  & $0$  & $0$  & $0$  & $0$ \\
		 $\chi_{34}$ & $0$  & $0$  & $1$  & $1$  & $0$  & $0$  & $0$  & -$1$  & -$1$  & $0$  & $0$  & $0$  & $0$  & $0$  & $0$  & $0$  & $0$  & $1$  & $1$  & -$1$  & -$1$  & -$1$  & $0$  & $0$ \\
		 $\chi_{35}$ & $C$  & $\overline{C}$  & $0$  & $0$  & $1$  & $0$  & $0$  & $0$  & $0$  & $0$  & $0$  & $0$  & $0$  & $0$  & $0$  & $C$  & $\overline{C}$  & $0$  & $0$  & $0$  & $0$  & $0$  & $0$  & $0$ \\
		 $\chi_{36}$ & $\overline{C}$  & $C$  & $0$  & $0$  & $1$  & $0$  & $0$  & $0$  & $0$  & $0$  & $0$  & $0$  & $0$  & $0$  & $0$  & $\overline{C}$  & $C$  & $0$  & $0$  & $0$  & $0$  & $0$  & $0$  & $0$ \\
		 $\chi_{37}$ & $0$  & $0$  & $3$  & $0$  & $0$  & $0$  & $0$  & $0$  & $0$  & $0$  & $0$  & $1$  & $1$  & $1$  & $0$  & $0$  & $0$  & $0$  & $0$  & -$1$  & -$1$  & -$1$  & -$1$  & -$1$ \\
		 $\chi_{38}$ & $2$  & $2$  & $0$  & $0$  & $0$  & $0$  & $0$  & $0$  & $0$  & $0$  & $0$  & -$1$  & -$1$  & -$1$  & $0$  & $0$  & $0$  & $1$  & $1$  & $0$  & $0$  & $0$  & -$1$  & -$1$ \\
		 $\chi_{39}$ & $0$  & $0$  & -$1$  & -$1$  & $0$  & $0$  & -$1$  & $0$  & $0$  & -$1$  & -$1$  & $0$  & $0$  & $0$  & $1$  & $0$  & $0$  & $0$  & $0$  & $1$  & $1$  & $1$  & $0$  & $0$ \\
		 $\chi_{40}$ & -$1$  & -$1$  & $2$  & -$1$  & $0$  & $0$  & $0$  & $0$  & $0$  & $0$  & $0$  & $0$  & $0$  & $0$  & $0$  & $1$  & $1$  & $0$  & $0$  & $0$  & $0$  & $0$  & -$1$  & -$1$ \\
		 $\chi_{41}$ & $0$  & $0$  & $0$  & $0$  & $0$  & $0$  & $0$  & $0$  & $0$  & $1$  & $1$  & $0$  & $0$  & $0$  & -$1$  & $0$  & $0$  & $0$  & $0$  & $0$  & $0$  & $0$  & $0$  & $0$ \\
		 $\chi_{42}$ & $0$  & $0$  & $0$  & $0$  & $0$  & $1$  & $0$  & $0$  & $0$  & $0$  & $0$  & $0$  & $0$  & $0$  & -$1$  & $0$  & $0$  & -$1$  & -$1$  & $0$  & $0$  & $0$  & $0$  & $0$ \\
		 $\chi_{43}$ & $0$  & $0$  & $0$  & $0$  & $1$  & $0$  & $1$  & -$1$  & -$1$  & $0$  & $0$  & $0$  & $0$  & $0$  & -$1$  & $0$  & $0$  & $0$  & $0$  & $0$  & $0$  & $0$  & $1$  & $1$ \\
		 $\chi_{44}$ & $0$  & $0$  & -$2$  & $1$  & $1$  & $0$  & -$1$  & $0$  & $0$  & $0$  & $0$  & $0$  & $0$  & $0$  & $0$  & $0$  & $0$  & $0$  & $0$  & $0$  & $0$  & $0$  & $0$  & $0$ \\
		 $\chi_{45}$ & $0$  & $0$  & $0$  & $0$  & $0$  & -$1$  & -$1$  & $1$  & $1$  & $0$  & $0$  & $0$  & $0$  & $0$  & $1$  & $0$  & $0$  & $0$  & $0$  & $0$  & $0$  & $0$  & -$1$  & -$1$ \\
		 $\chi_{46}$ & $0$  & $0$  & $0$  & $0$  & $0$  & $0$  & $0$  & $0$  & $0$  & $0$  & $0$  & $0$  & $0$  & $0$  & $0$  & $0$  & $0$  & $0$  & $0$  & $0$  & $0$  & $0$  & $0$  & $0$ \\
		 $\chi_{47}$ & $0$  & $0$  & $0$  & $0$  & $0$  & $0$  & $1$  & $0$  & $0$  & $0$  & $0$  & $0$  & $0$  & $0$  & -$1$  & $0$  & $0$  & $0$  & $0$  & $0$  & $0$  & $0$  & $1$  & $1$ \\
		 $\chi_{48}$ & $0$  & $0$  & $0$  & $0$  & -$1$  & $0$  & $0$  & $0$  & $0$  & $0$  & $0$  & $0$  & $0$  & $0$  & $1$  & $0$  & $0$  & $0$  & $0$  & $0$  & $0$  & $0$  & $0$  & $0$ \\
		 \bottomrule
	\end{tabular}
	\end{small}
	\end{center}
	\end{sidewaystable}

 \begin{sidewaystable}
\caption{Decompositions of $W^{(0)}_m, W^{(1)}_m$, part one. The representations appearing in the discriminant conjecture are in bold font.} \label{decomp1}
 \begin{center}
 \begin{small}
 \begin{tabular}{c@{ }@{ }r@{ }r@{ }r@{ }r@{ }r@{ }r@{ }r@{ }r@{ }r@{ }r@{ }r@{ }r@{ }r@{ }r@{ }r@{ }r@{ }r@{ }r@{ }r@{ }r@{ }r@{ }r@{ }r@{ }r}
 \toprule  $m \backslash$ & $V_{1}$& $V_{2}$& $V_{3}$&  ${\bf V_{4}}$&  ${\bf V_{5}}$& $V_{6}$& $V_{7}$& $V_{8}$& ${\bf V_{9}}$& ${\bf V_{10}}$&$ V_{11}$&$ V_{12}$& $V_{13}$& ${\bf V_{14}}$&  ${\bf V_{15}}$& $V_{16}$&  ${\bf V_{17}}$&  ${\bf V_{18}}$& $V_{19}$& $V_{20}$& $V_{21}$&  ${\bf V_{22}}$&  ${\bf V_{23}}$ & $V_{24}$ \\\midrule
 -3&-2&0&0&0&0&0&0&0&0&0&0&0&0&0&0&0&0&0&0&0&0&0&0&0 \\ 
 0& 0& 1& 0& 0& 0& 0& 0& 0& 0& 0& 0& 0& 0& 0& 0& 0& 0& 0& 0& 0& 0& 0& 0&0 \\
 1& 0& 0& 0& 0& 0& 0& 0& 0& 0& 0& 0& 0& 0& 0& 0& 0& 0& 0& 0& 0& 0& 0& 0&0 \\
 4& 0& 0& 0& 1& 1& 0& 0& 0& 0& 0& 0& 0& 0& 0& 0& 0& 0& 0& 0& 0& 0& 0& 0& 0 \\
 5& 0& 0& 0& 0& 0& 0& 0& 0& 1& 1& 0& 0& 0& 0& 0& 0& 0& 0& 0& 0& 0& 0& 0& 0 \\
 8& 0& 0& 0& 0& 0& 0& 0& 0& 0& 0& 0& 0& 0& 0& 0& 0& 1& 1& 0& 0& 0& 0& 0& 0 \\
 9& 0& 0& 0& 0& 0& 0& 0& 0& 0& 0& 0& 0& 0& 0& 0& 0& 0& 0& 0& 0& 0& 1& 1& 0 \\
 12& 0& 0& 0& 0& 0& 0& 0& 0& 0& 0& 0& 0& 0& 0& 0& 0& 0& 0& 0& 0& 0& 0& 0& 0 \\
 13& 0& 0& 0& 0& 0& 0& 0& 0& 0& 0& 0& 0& 0& 1& 1& 0& 0& 0& 0& 0& 0& 0& 0& 0 \\
 16& 0& 0& 0& 0& 0& 0& 0& 0& 0& 0& 0& 0& 0& 0& 0& 0& 0& 0& 0& 0& 0& 0& 0& 0 \\
 17 & 0& 0& 0& 0& 0& 0& 0& 0& 0& 0& 0& 0& 0& 0& 0& 0& 0& 0& 0& 0& 0& 0& 0& 0 \\
 20& 0& 0& 0& 0& 0& 0& 0& 0& 1& 1& 0& 0& 0& 0& 0& 0& 0& 0& 2& 0& 2& 0& 0& 0 \\
 21& 0& 0& 0& 0& 0& 0& 0& 0& 0& 0& 0& 0& 0& 0& 0& 0& 0& 0& 0& 0& 2& 2& 2& 2 \\
 24& 0& 0& 0& 0& 0& 0& 0& 0& 0& 0& 0& 0& 0& 2& 2& 4& 4& 4& 2& 6& 4& 8& 8& 8 \\
 25& 0& 0& 0& 0& 0& 0& 0& 0& 0& 0& 0& 2& 2& 2& 2& 2& 4& 4& 8& 8& 10& 12& 12& 12 \\
 28& 0& 0& 0& 2& 2& 0& 2& 2& 0& 0& 2& 14& 14& 8& 8& 10& 18& 18& 38& 28&  48& 54& 54& 56 \\
 29& 0& 0& 2& 0& 0& 0& 2& 2& 2& 2& 2& 18& 18& 18& 18& 20& 40& 40& 54& 58& 74& 86& 86& 86 \\
 32& 0& 0& 0& 0& 0& 4& 0& 4& 10& 10& 10& 62& 62& 78& 78& 92& 173& 173& 208& 256& 296& 368& 368& 368 \\
\bottomrule
\end{tabular}
\end{small}
\end{center}
\end{sidewaystable}
	
 \begin{sidewaystable}
\caption{Decompositions of $W^{(0)}_m, W^{(1)}_m$, part two.} \label{decomp2}
 \begin{center}
 \begin{small}
 \begin{tabular}{c@{ }@{ }r@{ }r@{ }r@{ }r@{ }r@{ }r@{ }r@{ }r@{ }r@{ }r@{ }r@{ }r@{ }r@{ }r@{ }r@{ }r@{ }r@{ }r@{ }r@{ }r@{ }r@{ }r@{ }r@{ }r}
 \toprule $m \backslash$ & $V_{25}$& $V_{26}$&  ${\bf V_{27}}$&  ${\bf V_{28}}$&  ${\bf V_{29}}$&  ${\bf V_{30}}$& $V_{31}$& $V_{32}$& $V_{33}$&$V_{34}$& ${\bf V_{35}}$& ${\bf V_{36}}$& $V_{37}$ &$V_{38}$& $V_{39}$& $V_{40}$& $V_{41}$& $V_{42}$& $V_{43}$& $V_{44}$& $V_{45}$& $V_{46}$& $V_{47}$& $V_{48}$ \\\midrule
 -3&0&0&0&0&0&0&0&0&0&0&0&0&0&0&0&0&0&0&0&0&0&0&0&0 \\ 
 0& 0& 0& 0& 0& 0& 0& 0& 0& 0& 0& 0& 0& 0& 0& 0& 0& 0& 0& 0& 0& 0& 0& 0&0 \\
 1& 0& 0& 0& 0& 0& 0& 0& 0& 0& 0& 0& 0& 0& 0& 0& 0& 0& 0& 0& 0& 0& 0& 0&0 \\
 4& 0& 0& 0& 0& 0& 0& 0& 0& 0& 0& 0& 0& 0& 0& 0& 0& 0& 0& 0& 0& 0& 0& 0&0 \\
 5& 0& 0& 0& 0& 0& 0& 0& 0& 0& 0& 0& 0& 0& 0& 0& 0& 0& 0& 0& 0& 0& 0& 0&0 \\
 8& 0& 0& 0& 0& 0& 0& 0& 0& 0& 0& 0& 0& 0& 0& 0& 0& 0& 0& 0& 0& 0& 0& 0&0 \\
 9& 0& 0& 0& 0& 0& 0& 0& 0& 0& 0& 0& 0& 0& 0& 0& 0& 0& 0& 0& 0& 0& 0& 0&0 \\
 12& 0& 0& 0& 0& 0& 0& 0& 0& 0& 0& 0& 0& 0& 0& 0& 2& 0& 0& 0& 0& 0& 0& 0& 0 \\
 13& 0& 0& 0& 0& 0& 0& 0& 0& 0& 0& 0& 0& 0& 0& 0& 0& 0& 0& 0& 0& 0& 2& 0& 0 \\
 16 & 0& 0& 1& 1& 0& 0& 0& 0& 0& 0& 0& 0& 0& 0& 0& 0& 0& 2& 2& 0& 2& 2& 0& 4 \\
 17 & 0& 0& 0& 0& 0& 0& 0& 0& 0& 0& 0& 0& 2& 2& 2& 0& 2& 2& 2& 2& 2& 2& 4& 6 \\
 20& 2& 0& 0& 0& 0& 0& 0& 2& 4& 4& 4& 4& 8& 6& 12& 10& 12& 10& 12& 16& 14& 16& 24& 32 \\
 21& 2& 2& 2& 2& 2& 2& 4& 6& 6& 8& 8& 8& 8& 10& 12& 16& 14& 26& 24& 26& 26& 28& 36& 62 \\
 24& 4& 10& 14& 14& 14& 14& 22& 18& 22& 28& 38& 38& 44& 52& 56& 68& 80& 132& 136& 126& 140& 154& 180& 332 \\
 25& 14& 14& 20& 20& 20& 20& 32& 34& 50& 54& 60& 60& 82& 86& 116& 128& 150& 210& 220& 222& 234& 264& 320& 550 \\
 28& 76& 64& 82& 82& 84& 84& 140& 160& 242& 256& 284& 284& 388& 404& 560&  616& 706& 950& 998& 1048& 1076& 1220& 1504& 2510 \\
 29& 104& 106& 144& 144& 144& 144& 232& 248& 352& 390& 462& 462& 636& 668& 880& 960& 1130& 1582& 1664& 1684& 1768& 1986& 2428& 4146 \\
 32& 416& 450& 624& 624& 624& 624& 992& 1040& 1448& 1616& 1932& 1932& 2620& 2780& 3602& 3956& 4666& 6686& 6978& 6988& 7394& 8278& 10068& 17372 \\
\bottomrule
\end{tabular}
\end{small}
\end{center}
\end{sidewaystable}
\vfil \eject

\end{document}